\documentclass[12pt]{amsart}
\usepackage{amscd, amsfonts, amssymb, amsthm}
\usepackage[all]{xypic}
\usepackage{graphics, lscape}
\usepackage{rotating}

\newcommand\CB{{\mathcal B}} 

\renewcommand\CD{{\mathcal D}} 

\newcommand\CH{{\mathcal H}} 
\newcommand\CL{{\mathcal L}}
 
\newcommand\CO{{\mathcal O}}

\newcommand\CS{{\mathcal S}}
\newcommand\CT{{\mathcal T}}

\newcommand\BBC{{\mathbb C}}

\newcommand\BBF{{\mathbb F}}
\newcommand\BBN{{\mathbb N}}

\newcommand\BBZ{{\mathbb Z}}

\newcommand\bs{{\mathbf s}}
\newcommand\bx{{\mathbf x}}

\newcommand\Atilde{{\widetilde{A}}}

\newcommand\CHtilde{{\widetilde{\CH}}}

\newcommand\End{{\operatorname{End}}}
\newcommand\GL{{\operatorname{GL}}}

\newcommand\ind{{\operatorname{ind}}}

\newcommand\inverse{^{-1}}

\renewcommand\a{\mathbf a}
\newcommand\abar{\overline{\a}}
\newcommand\CHt{\widehat{\mathcal H}}
\newcommand\CHqa{\mathcal H_{q,a}}
\newcommand\CHv{\mathcal H_{\mathbf v}}
\newcommand\q{\mathbf q}
\newcommand\Tsbar{\overline{T}_s}
\newcommand\bv{\mathbf v}
\newcommand\vbar{\overline{\mathbf v}}

\makeatletter
\newdimen\p@renwd
\setbox0=\hbox{\rm B} \p@renwd=\wd0 
\def\bordermatrix#1{\begingroup \m@th
  \setbox\z@\vbox{\def\cr{\crcr\noalign{\kern2\p@\global\let\cr\endline}}%
    \ialign{$##$\hfil\kern2\p@\kern\p@renwd&\thinspace\hfil$##$\hfil
      &&\quad\hfil$##$\hfil\crcr
      \omit\strut\hfil\crcr\noalign{\kern-\baselineskip}%
      #1\crcr\omit\strut\cr}}%
  \setbox\tw@\vbox{\unvcopy\z@\global\setbox\@ne\lastbox}%
  \setbox\tw@\hbox{\unhbox\@ne\unskip\global\setbox\@ne\lastbox}%
  \setbox\tw@\hbox{$\kern\wd\@ne\kern-\p@renwd\left[\kern-\wd\@ne
    \global\setbox\@ne\vbox{\box\@ne\kern2\p@}%
    \vcenter{\kern-\ht\@ne\unvbox\z@\kern-\baselineskip}\,\right]$}%
  \null\;\vbox{\kern\ht\@ne\box\tw@}\endgroup}
\makeatother

\numberwithin{equation}{subsection}

\theoremstyle{plain}
\newtheorem{lemma}[subsection]{Lemma}
\newtheorem{theorem}[subsection]{Theorem}
\newtheorem{corollary}[subsection]{Corollary}
\newtheorem{proposition}[subsection]{Proposition}

\theoremstyle{definition}
\newtheorem{example}[subsection]{Example}

\newtheorem*{lemma*}{Lemma}
\newtheorem*{proposition*}{Proposition}
\newtheorem*{theorem*}{Theorem}
\newtheorem*{notation*}{Notation}

\begin{document}


\title[Algebras for Monomial Groups] {Generic Hecke Algebras for
  Monomial Groups}

\author[S.I. Alhaddad]{S.I. Alhaddad} \address{Department of
  Mathematics\\ University of South Carolina, Lancaster\\ Lancaster,
  SC 29721}
\email{alhaddad@gwm.sc.edu}

\author[J.M. Douglass]{J. Matthew Douglass} \address{Department of
  Mathematics\\ PO Box 311430\\ University of North Texas\\ Denton TX,
  USA 76203-1430} \email{douglass@unt.edu}
 


\begin{abstract}
  In this paper we define a two-variable, generic Hecke algebra, $\mathcal
  H$, for each complex reflection group $G(b,1,n)$. The algebra $\mathcal H$
  specializes to the group algebra of $G(b,1,n)$ and also to an endomorphism
  algebra of a representation of $\operatorname{GL}_n(\mathbb F_q)$ induced
  from a solvable subgroup. We construct Kazhdan-Lusztig ``$R$-polynomials''
  for $\CH$ and show that they may be used to define a partial order on
  $G(b,1,n)$. Using a generalization of Deodhar's notion of distinguished
  subexpressions we give a closed formula for the $R$-polynomials. After
  passing to a one-variable quotient of the ring of scalars, we construct
  Kazhdan-Lusztig polynomials for $\mathcal H$ that reduce to the usual
  Kazhdan-Lusztig polynomials for the symmetric group when $b=1$.
\end{abstract}

\thanks{The authors would like to thank Nathaniel Thiem for helpful
  discussions.}

\maketitle
\allowdisplaybreaks


\section{Introduction}\label{section1}

Suppose $G$ is a finite group of Lie type and $k$ is a field whose
characteristic is relatively prime to $|G|$. Then the irreducible
$k$-representations of $G$ are partitioned into Harish-Chandra series (see
for example \cite{curtisreiner:methodsII}). An irreducible representation is
cuspidal if it is the only member in its Harish-Chandra series. Non-cuspidal
representations can be constructed and their character values computed using
the theory of endomorphism algebras of induced representations.

The prototypical example of a Harish-Chandra series that does not consist of
just one cuspidal representation is the principal series.  In this case, the
endomorphism algebra of interest is the endomorphism algebra of a
representation induced from the trivial representation of a Borel subgroup
of $G$. This algebra is a specialization of a generic algebra, the
Iwahori-Hecke algebra of the Weyl group $W$ of $G$. The Iwahori-Hecke
algebra also specializes to the group algebra of $W$.

Brou\'e, Malle, and Michel (see \cite{brouemalle:zyklotomische},
\cite{brouemallemichel:generic}, \cite{brouemichel:blocs}) have shown that
if the characteristic of $k$ is different from the defining characteristic
of $G$, then the classical Harish-Chandra theory just outlined can be
extended to describe the blocks of $kG$. A new feature that arises when the
characteristic of $k$ divides $|G|$ is that the endomorphism rings used to
analyze non-cuspidal representations are deformations of group algebras of
complex reflection groups that are not Weyl groups or even Coxeter groups.

In order to study Iwahori-Hecke algebras, Kazhdan and Lusztig
\cite{kazhdanlusztig:coxeter} and Lusztig \cite{lusztig:characters} have
developed a powerful theory for analyzing representations of Iwahori-Hecke
algebras which in turn plays a central role in describing the irreducible
representations of $G$.

In this paper we consider a variation of the above themes and define a
two-variable generic algebra, $\CH$, for finite general linear groups that
specializes to the group algebra of the complex reflection group $G(b,1,n)$
and also specializes to the endomorphism algebra of an induced
representation of $\GL_n(\BBF_q)$. This last algebra is closely related to
the principal block of $\GL_n(\BBF_q)$ when $b=|\GL_n(\BBF_q)|_l$ where $l$
is the characteristic of $k$ (see \S3).

Starting from the observation that the Bruhat-Chevalley order on a finite
Coxeter group is determined by the non-vanishing of the ``$R$-polynomials''
of Kazhdan and Lusztig, we construct $R$-polynomials for $\CH$ and show that
the resulting relation on $G(b,1,n)$ is a partial order. This partial order
is the Bruhat order on the symmetric group when $b=1$ but it is not the
Bruhat order on hyperoctahedral groups when $b=2$.

Passing to a single variable quotient of our two-variable ring of scalars we
define Kazhdan-Lusztig polynomials for $G(b,1,n)$. When $b=1$ these are the
usual Kazhdan-Lusztig polynomials for the symmetric group.

Our construction has a generalization to arbitrary finite Chevalley groups
or $\BBF_q$-points of a reductive algebraic group. When the underlying root
system is not of type $A$, the groups that arise are no longer complex
reflection groups. In order to retain the connection with complex reflection
groups and to make the exposition of the ideas as accessible as possible, in
this paper we consider only the case of the general linear group.

A specialization of the algebra $\CH$ has been considered by Cabanes and
Enguehard \cite[Chapter 23] {cabanesenguehard:representations}. Because our
construction is generic it opens up the possibility of relating the
specialized algebra studied by Cabanes and Enguehard with other algebras
used to study representations of finite reductive groups. In particular,
because $\CH$ specializes to the group algebra of $G(b,1,n)$ it should be
related to a suitable specialization of an Ariki-Koike algebra. It also
seems likely that $\CH$ is a cellular algebra in the sense of Graham and
Lehrer \cite{grahamlehrer:cellular}. We hope to consider these questions, as
well as the representation theory of $\CH$, in future work.

For the rest of this paper, $G$ will denote $\GL_n(\BBF_q)$ where $q$ is a
prime power. Suppose that $a$ and $b$ are relatively prime, positive
integers with $ab=q-1$. The multiplicative group of $\BBF_q$ is a cyclic
group with order $q-1$ and so it factors as the direct product of a cyclic
group of order $a$, which we will denote by $F_a$, and a cyclic group of
order $b$, which we will denote by $F_b$. Let $H$ be the subgroup of $G$
consisting of diagonal matrices. Let $H_a$ and $H_b$ be the subgroups of $H$
with entries in $F_a$ and $F_b$ respectively. Then clearly $H\cong H_a\times
H_b$.

Let $U$ denote the subgroup of $G$ consisting of upper triangular,
unipotent matrices and define $B_a=H_aB$. Notice that $B_a=U$ when
$a=1$ and $B_a$ is a Borel subgroup of $G$ when $a=q-1$. 

Let $e$ denote the centrally primitive idempotent in the group algebra
$\BBC B_a$ corresponding to the trivial representation of $B_a$.  In
\S2 we study the subalgebra $e\BBC Ge$. This algebra is isomorphic to
the opposite algebra of the endomorphism algebra of the representation
of $G$ induced from the trivial representation of $B_a$. Let $W$
denote the subgroup of $G$ consisting of permutation matrices. It is
easily seen that the subgroup $WH_b$ of $G$ is isomorphic to
$G(b,1,n)$ and that $WH_b$ is in fact a complete set of $(B_a,
B_a)$-double coset representatives.  In \S2 we prove some
multiplication relations in $e\BBC Ge$.  These relations are analogous
to the braid and quadratic relations in the Iwahori-Hecke algebra of
$W$ and also to the relations in the Bernstein-Zelevinsky presentation
of the extended, affine Hecke algebra of $W$.

In \S3 we define a generic algebra $\CH$ using the relations from \S2 as a
model. This is entirely analogous to the construction of the Iwahori-Hecke
algebra. We show that our generic algebra has a basis indexed by $WH_b$
using an argument that goes back to \cite[Ch.~4,
Ex.~23]{bourbaki:groupes}. We then record some of the standard properties of
the Iwahori-Hecke algebra that remain true in our setup and give an
application to the characters in the principal block of $G$ in the case when
$H_b$ is a Sylow-$l$ subgroup of $G$ and $l$ is the characteristic of $k$.

In the general case when $G_0$ is a finite Chevalley group with defining
characteristic equal $p$ and $U$ is a Sylow-$p$ subgroup of $G_0$, Yokonuma
\cite{yokonuma:structure} has given a presentation by generators and
relations for the Hecke ring $\CH(G_0,U)$ and defined such generic algebra
for $(G_0,U)$. The Hecke ring $\CH(G,P)$, for a subgroup $P$ of $G_0$
containing $U$, is then realized as a subalgebra of the Hecke ring
$\CH(G_0,U)$ obtained by summing basis elements of $\CH(G_0,U)$. Our results
in \S2 and \S3 may alternatively be derived from Yokonuma's results in the
same way that the Hecke ring $\CH(G,B)$ and the Iwahori-Hecke algebra are
derived from the group algebra of $G_0$ by analyzing the structure constants
for a particular choice of basis.

In \S4 we construct Kazhdan-Lusztig ``$R$-polynomials'' and show that they
can be used to define a partial order on $WH_b$. In the proof, we relate
intervals in the partial order with intervals in the Bruhat order on $W$ and
certain subsets of $H_b$.

The algebra $\CH$ is a $\BBZ[\a, \bv,\bv\inverse]$-algebra where $\a$ and
$\bv$ are indeterminates. In \S5, after passing to a quotient isomorphic to
$\BBZ_b[\bv, \bv\inverse]$, where $\BBZ_b$ denotes the localization of
$\BBZ$ at $b$, we construct a Kazhdan-Lusztig basis and Kazhdan-Lusztig
polynomials following an argument in \cite{lusztig:hecke}. We compute these
``Kazhdan-Lusztig'' polynomials when $n=3$ and $b$ is arbitrary. When $b=1$,
these are the usual Kazhdan-Lusztig polynomials and so are known a priori to
be the constant polynomial $1$ in all cases. In contrast, when $b\ne1$, some
of our polynomials have positive degree. 

Finally, in \S6 we adapt Deodhar's ideas in \cite{deodhar:parabolic}
to describe the lower order ideals in the poset $(WH_b, \leq)$ and to
give a closed form expression for the $R$-polynomials. These arguments
use an analog in $WH_b$ of a reduced expression of an element in $W$.

\section{The Algebra $e\BBC Ge$}\label{section2}

Recall that $G=\GL_n(\BBF_q)$ and that $q-1=ab$ where $a$ and $b$ are
relatively prime. 

In this section $\CHqa$ will denote the ``Hecke algebra'' $e\BBC Ge$,
where $\BBC G$ is the group algebra of $G$ and $e$ is the centrally
primitive idempotent $|B_a|\inverse \sum_{b\in B_a} b$ in $\BBC B_a$.
Thus, $\CH_{q,a}$ is anti-isomorphic to the endomorphism ring of the
induced representation $\operatorname {Ind}_{B_a}^G( 1_{B_a})$.

Suppose $\{\, x_1, \dots, x_m\,\}$ is a complete set of $(B_a, B_a)$
double coset representatives and $D_i= B_ax_iB_a$ is the double coset
containing $x_i$. It is well known (see \cite[Proposition
11.34]{curtisreinermethodsI}) that if we consider $\BBC G$ as
$\BBC$-valued functions on $G$ and let $T_{x_i}= |B_a|\inverse
\chi_i$, where $\chi_i$ is the characteristic function of $D_i$, then
the $T_{x_i}$'s are a basis of $\CHqa$.  Moreover, the multiplication
in $\CHqa$ is given by
\[
T_{x_i} T_{x_j}= \sum_{k=1}^m \mu_{x_i,x_j,x_k} T_{x_k},\ \text{where
  $\mu_{x_i,x_j,x_k}= |B_a|\inverse |D_i \cap x_k D_j\inverse|$.}
\]

For a permutation matrix $w$ in $W$, define $U_{w}^-= \{\,u\in U\mid
wuw\inverse \in w_0Uw_0\,\}$, where $w_0$ is the permutation matrix
with $1$'s on the antidiagonal. Then by the strong form of the Bruhat
decomposition for $G$ we have $G=\coprod_{w\in W} U_{w\inverse}^-wHU$
with uniqueness of expression. Since $H\cong H_b \times H_a$, this
proves the next lemma.

\begin{lemma}
  Every element in $G$ has a unique expression as a product
  $u_1wt_bt_au_2$ where $w$ is in $W$, $u_1$ is in $U_{w\inverse}^-$,
  $h_a$ is in $H_a$, $h_b$ is in $H_b$, and $u$ is in $U$. In
  particular, the subgroup $WH_b$ of $G$ is a complete set of $(B_a,
  B_a)$-double coset representatives.
\end{lemma}

It follows that $\dim \CHqa= |WH_b|= b^n n!$. 

Fix a generator, $\zeta$, of $\BBF_q^\times$. Then $\zeta^b$ generates
$F_a$ and $\zeta^a$ generates $F_b$.  In order to determine the
structure constants $\mu_{x,y,z}$ for $x$, $y$, and $z$ in $WH_b$ we
need the following lemma, whose easy proof will be omitted.

\begin{lemma}
  Suppose $u= \left[\begin{smallmatrix} 1&\zeta^c\\ 0&1
      \end{smallmatrix} \right]$ is in $\GL_2(\BBF_q)$ with $0<c<q-1$
    and $s= \left[\begin{smallmatrix} 0&1\\ 1&0 \end{smallmatrix}
    \right]$. Write $c=am+bn$ where $0\leq m\leq b-1$ and $0\leq n\leq
    a-1$.  Then $sus= u_1sh_bh_au_2$ where
  \[
  u_1=u_2= \begin{bmatrix} 1&\zeta^{-c}\\ 0&1 \end{bmatrix}, \quad
  h_b= \begin{bmatrix} \zeta^{am}&0\\ 0&(-1)^{b-1}\zeta^{-am}
  \end{bmatrix}, \quad\text{and}\quad h_a= \begin{bmatrix}
    \zeta^{bn}&0\\ 0&(-1)^b\zeta^{-bn} \end{bmatrix}. 
  \]
\end{lemma}

Notice that $-1= (\zeta^b)^{a/2}$ is in $F_a$ if $b$ is odd and $-1=
(\zeta^a)^{b/2}$ is in $F_b$ if $b$ is even.

For $1\leq i\leq n-1$  define 
\[
s_i=\bordermatrix{&&& i & i+1 && \cr &1&&&&& \cr && \ddots& & & & \cr
  i&&&0 &1 && \cr i+1&&&1 & 0 && \cr &&& & & \ddots &\cr &&&&&&&1}
\]
so $\CS= \{\, s_1, \dots, s_{n-1}\,\}$ is a set of Coxeter generators
for $W$.

For $\alpha$ in $F_b$ and $1\leq i\ne j\leq n$, define $h_i(\alpha)$
to be the diagonal matrix whose $i$th entry is $\alpha$ and whose
other diagonal entries are $1$ and define $h_{i,j} (\alpha)=
h_i(\alpha) h_j((-1)^{b-1}\alpha \inverse)$.

\begin{corollary}\label{sl2}
  Suppose $s=s_i$ is in $S$.
  \begin{enumerate}
  \item If $u$ is in $U_s^-$ and the non-zero, off-diagonal entry of
    $u$ is $\zeta^c$ where $c=am+bn$, then $sus= u_1sh_b^uh_a^uu_2$
    where $h_b^u=h_{i, i+1}(\zeta^{am})$.
  \item The relation $u_1\sim u_2$ if and only if $h_b^{u_1}=
    h_b^{u_2}$ is an equivalence relation on $U_s^-$ and each
    equivalence class contains $a$ elements.
  \end{enumerate} 
\end{corollary}

\begin{proof}
  The first statement follows immediately from the lemma. 
  
  The relation in the second statement is clearly an equivalence
  relation. If $u$ is in $U_s^-$ and the non-zero, off-diagonal entry
  of $u$ is $\zeta^{am+bn}$, then it follows from the first statement
  that the equivalence class containing $u$ is the set of all elements
  in $U_s^-$ whose non-zero, off-diagonal entry is $\zeta^{am+bn'}$
  where $0\leq n'\leq a-1$.
\end{proof}

Using the natural projection $WH_b \to W$ we can lift the length
function, $\ell\colon W\to \BBN$, to $WH_b$. Then $\ell(wd)=
\ell(dw)=\ell(w)$ for $w$ in $W$ and $d$ in $H_b$.

Let $\CT= \{\, ws_iw\inverse \mid 1\leq i\leq n-1,\, w\in W\,\}$ be
the set of reflections in $W$.  For $t=ws_iw\inverse$ in $\CT$, define
\begin{equation}
  \label{xs}
  X_t=\{\, wh_{i, i+1}(\alpha) w\inverse \mid \alpha\in F_b \,\}.
\end{equation}
Then $X_t$ does not depend on the choice of $(w,i)$. Also, if $b$ is
odd, then $X_t$ is a subgroup of $H_b$ isomorphic to $F_b$; if $b$ is
even, then $X_t$ is a coset of a subgroup of $H_b$ isomorphic to
$F_b$.

We can now describe the multiplication in $\CHqa$.

\begin{theorem}\label{multH}
  The multiplication in the Hecke algebra $\CHqa$ is determined by the
  following relations.
  \begin{enumerate}
  \item If $d$ and $d'$ are in $H_b$, then $T_{d'} T_{d}= T_{d'd}$.
    Thus, the rule $d\mapsto T_d$ extends linearly to an algebra
    isomorphism between the group algebra $\BBC H_b$ and the subset of
    $\CHqa$ spanned by $\{\, T_d\mid d\in H_b\,\}$.
  \item If $d$ is in $H_b$ and $w$ is in $W$, then $T_wT_d= T_{wd}$
    and $T_dT_w= T_{dw}$.
  \item If $w$ is in $W$ and $s$ is in $\CS$, then
    \[
    T_w T_s=\begin{cases} T_{ws}& \text{if $\ell(ws)> \ell(w)$} \\
      qT_{ws}+ a\sum_{d\in X_s} T_{wd} & \text{if $\ell(ws)<
        \ell(w)$.} 
    \end{cases} 
    \]
  \end{enumerate}
\end{theorem}

\begin{proof}
  Recall that for $x$ and $y$ in $WH_b$ we have 
  \[
  T_xT_y= \sum_{z\in WH_b} \mu_{x,y,z}T_z 
  \]
  where $\mu_{x,y,z}= |B_a|\inverse |B_axB_a \cap zB_ay\inverse B_a|$. The
  statements in the theorem are all proved by computing $\mu_{x,y,z}$ for
  appropriate choices of $x$, $y$ and $z$. The most complicated case is
  showing that $T_w T_s= qT_{ws}+ a\sum_{d\in X_s} T_{wd}$ when $\ell(ws)<
  \ell(w)$. We prove this statement and omit the other computations.

  Fix $w$ in $W$ and $s=s_i$ in $S$ with $\ell(ws)< \ell(w)$. Then for $x$
  in $WH_b$ we have
  \begin{align*}
    \mu_{w,s,x}&=\frac 1 {|B_a|} \left| B_awB_a \cap xB_asB_a\right|\\
    &= \frac 1 {|B_a|} \left| \{\, (u,b)\in U_s^- \times B_a \mid
      xusb\in B_awB_a \,\} \right|\\
    &= \left| \{\, u\in U_s^- \mid xus \in B_awB_a \,\}
    \right|.
  \end{align*}

  Suppose $x=w_1d$ where $w_1$ is in $W$ and $d$ is in $H_b$.

  If $\ell(w_1s)> \ell(w_1)$, then for $u$ is $U_s^-$ we have $w_1dus=
  u'w_1sd'$ where $u'=w_1du d\inverse w_1\inverse$ is in $U$ and $d'=sds$ is
  in $H_b$. Thus, $xus$ is in $B_awB_a$ if and only if $w_1=ws$ and
  $d=1$. Therefore $\mu_{w,s,ws}=q$ and $\mu_{w,s,x}=0$ otherwise.

  If $\ell(w_1s)< \ell(w_1)$, then $w_1dus= w_1s (sds)$ is in $B_awB_a$ if
  and only if $w_1=ws$ and $d=1$. But then $\ell(w_1s)= \ell(w)> \ell(ws)=
  \ell(w_1)$, a contradiction, so $\mu(w,s,x)=0$ in this case.  If
  $\ell(w_1s)< \ell(w_1)$ and $u$ is in $U_s^-$ with $u\ne 1$, then $w_1dus=
  u'w_1dh_b^ub'$ where $u'=w_1dsu_1sd\inverse w_1\inverse$ is in $U$,
  $b'=h_a^uu_2$ is in $B_a$, and $sus=u_1s h_b^uh_a^uu_2$ as in Corollary
  \ref{sl2}. Thus, $xus$ is in $B_awB_a$ if and only if $w_1=w$ and
  $d=(h_b^u)\inverse$. It follows that $\mu_{w,s,x}=0$ unless $w_1=w$ and
  $d$ is in $X_s$. For $x=wd$ with $d$ in $X_s$ we have
  \[ 
  \mu_{w,s,wd}= |\{\, u\in U_s^-\mid d=(h_b^u)\inverse\,\}|.
  \]
  It follows from Corollary \ref{sl2} that for a given $h_{i, i+1}
  (\alpha)$, there are $a$ elements, $u$, in $U_s^-$ with $h_b^u=
  h_{i, i+1} (\alpha)$ and so $\mu_{w,s,wd}=a$ for $d$ in $X_s$.

  We have shown that if $\ell(ws)<\ell(w)$, then
  \[
  \mu_{w,s,x}=
  \begin{cases}
    q& x= ws\\
    a&x=wd,\, d\in X_s\\
    0&\text{otherwise}
  \end{cases}
  \]
  and so $T_wT_s= qT_{ws}+ a\sum_{d\in X_s} T_{wd}$. This completes
  the proof of the theorem.
\end{proof}

We conclude this section by recording some properties of $H_b$, and the
subgroups $X_t$ will be used later.

\begin{proposition}\label{Xt}
  For every $t$ in $\CT$, the subset $X_t$ is closed under taking
  inverses. If $t_1$, \dots, $t_{r+1}$ are in $\CT$, then
  \[
  X_{t_1} \dotsm X_{t_r} X_{t_{r+1}} = X_{t_1} \dotsm X_{t_r}
  X_{t_1\dotsm t_rt_{r+1} t_r \dotsm t_1}.
  \]
\end{proposition}

\begin{proof}
  The fact that $X_t=X_t\inverse$ follows immediately from the
  definition. 

  To prove the second statement we use induction on $r$.
  
  Suppose $r=1$. If $t_1t_2=t_2t_1$, then the result is clear. Suppose
  $t_1t_2\ne t_2t_1$. There are distinct $i$, $j$, and $k$ with $1\leq
  i,j,k\leq n$ so that $t_1$ interchanges the $i$th and $j$th
  standard basis vectors of $\BBF_q^n$, $t_2$ interchanges the $j$th
  and $k$th standard basis vectors of $\BBF_q^n$, and $t_1t_2t_1$
  interchanges the $i$th and $k$th standard basis vectors of
  $\BBF_q^n$. Then $h_{i,j}(\alpha) h_{j,k}(\beta) =h_{i,j}(\alpha
  \beta\inverse) h_{i,k}(\beta)$ for $\alpha$ and $\beta$ in $F_b$ and
  so
  \begin{align*}
    X_{t_1}X_{t_2}&= \{\, h_{i,j}(\alpha) h_{j,k}(\beta) \mid \alpha,
    \beta \in F_b \,\}\\
    &= \{\, h_{i,j}(\alpha\beta\inverse) h_{i,k}(\beta) \mid \alpha,
    \beta \in F_b \,\}\\
    &= \{\, h_{i,j}(\alpha_1) h_{i,k}(\beta_1) \mid \alpha_1,
    \beta_1 \in F_b \,\}\\
    &= X_{t_1}X_{t_1t_2t_1}.
  \end{align*}
  
  Now suppose that $r>1$. Then by induction and the case when $r=1$ we
  have
  \begin{align*}
    X_{t_1} X_{t_2} \dotsm X_{t_r} X_{t_{r+1}} &= X_{t_1}
    X_{t_2}\dotsm
    X_{t_r} X_{t_2\dotsm t_rt_{r+1} t_r \dotsm t_2}\\
    &= X_{t_2}\dotsm X_{t_r} X_{t_1} X_{t_2\dotsm t_rt_{r+1} t_r
      \dotsm
      t_2}\\
    &= X_{t_2}\dotsm X_{t_r} X_{t_1} X_{t_1\dotsm t_rt_{r+1} t_r
      \dotsm
      t_1}\\
    &=X_{t_1} \dotsm X_{t_r} X_{t_1\dotsm t_rt_{r+1} t_r \dotsm t_1}.
  \end{align*}
  This completes the proof of the proposition.
\end{proof}

Define $X_0=\{\, h_1(\zeta^{ai}) \mid 0\leq i\leq b-1\,\}$, so $X_0$
is a subgroup of $H_b$.

\begin{proposition}\label{product}
  The multiplication mapping
  \[
  X_0\times X_{s_1} \times \dotsm \times X_{s_{n-1}} \to H_b
  \]
  is a bijection.
\end{proposition}

\begin{proof}
  Since the domain and codomain of the mapping both have cardinality
  $b^n$, it is enough to show that the mapping is injective. 

  We use induction on $n$, the base case being when $n=2$.

  Suppose $h_1(\zeta^{ai}) h_{1,2}(\alpha_1) = h_1(\zeta^{aj})
  h_{1,2}(\beta_1)$. Then comparing $(2,2)$-entries we see that
  $(-1)^{b-1} \alpha_{1} = (-1)^{b-1} \beta_{1}$ and so
  $\alpha_{1}= \beta_{1}$. It follows that $h_1(\zeta^{ai})=
  h_1(\zeta^{aj})$.

  Suppose $n>2$ and  
  \[
  h_1(\zeta^{ai}) h_{1,2}(\alpha_1) \dotsm h_{n-1,n}(\alpha_{n-1}) =
  h_1(\zeta^{aj}) h_{1,2}(\beta_1) \dotsm h_{n-1,n}(\beta_{n-1}).  
  \]
  Then comparing $(n,n)$-entries we see that $(-1)^{b-1} \alpha_{n-1}
  = (-1)^{b-1} \beta_{n-1}$ and so $\alpha_{n-1}= \beta_{n-1}$.
  Therefore
  \[
  h_1(\zeta^{ai}) h_{1,2}(\alpha_1) \dotsm h_{n-2,n-1}(\alpha_{n-2}) =
  h_1(\zeta^{aj}) h_{1,2}(\beta_1) \dotsm h_{n-2,n-1}(\beta_{n-2}).  
  \]
  By induction we have $h_1(\zeta^{ai})= h_1(\zeta^{aj})$ and
  $\alpha_i=\beta_i$ for $1\leq i\leq n-2$.
\end{proof}

\section{The Generic Algebra $\CH$}

In this section we define a generic algebra that specializes to the
algebra $\CHqa$ from \S\ref{section2}. This algebra depends only on
the triple $(W,\CS,b)$ where $(W,\CS)$ is a Coxeter group with $W$ a
symmetric group, and $b$ is a positive integer. Such a triple
determines the wreath product $(\BBZ/b\BBZ)\wr W$. In order to keep
the notation to a minimum we will always consider the particular
representation of this group as the subgroup $WH_b$ of $\GL_n(\BBF_q)$
from \S\ref{section2}. In addition, we will continue to use the
notation already introduced for $WH_b$. In particular, $\ell$ is the
length function, $\CS$ is a fixed set of Coxeter generators of $W$,
and for $s$ in $\CS$, $X_s$ is the subset of $H_b$ defined in
(\ref{xs}).

Let $A=\BBZ[\a,\bv]$ where $\a$ and $\bv$ are indeterminates. Set
$\q=\bv^2$ and define $\CH$ to be the $A$-algebra with generators $\{\,
t_s\mid s\in \CS\,\} \cup \{\, t_d\mid d\in H_b\,\}$ and relations
\begin{gather}
  t_dt_{d'}= t_{dd'}\label{rel1}\\
  t_dt_{s_i}=t_{s_i}t_{s_ids_i}\label{rel1.5}\\
  t_{s_i}t_{s_j}=t_{s_j}t_{s_i}\quad \text{if $|j-i|>1$}\label{rel2}\\
  t_{s_i}t_{s_{i+1}}t_{s_i}= t_{s_{i+1}}t_{s_i}t_{s_{i+1}}
  \label{rel3}\\ 
  t_{s_i}^2=\q1_{\CH}+ \a \sum_{d\in X_{s_i}} t_dt_{s_i} \label{rel4}
\end{gather}
where $1\leq i<j\leq n-1$ and $d$ and $d'$ are in $H_b$. 

It follows from relations (\ref{rel2}) and (\ref{rel3}) and
Matsumoto's Theorem \cite[theorem 1.2.2]{geckpfeiffer:characters} that
if $w$ is in $W$, $w=s_{i_1} \dotsm s_{i_p}$, and $\ell(w)=p$, then
the product $t_{s_1} \dotsm t_{s_p}$ depends only on $w$ and not on
$s_{i_1}$, \dots, $s_{i_p}$ so we may define $t_w$ unambiguously by
$t_w=t_{s_1} \dotsm t_{s_p}$. We also define $t_{wd}=t_wt_d$ and
$t_{dw}=t_dt_w$ for $w$ in $W$ and $d$ in $H_b$. Using relation
(\ref{rel1.5}) and induction it is easy to see that $t_x$ is
unambiguously defined for every $x$ in $WH_b$. 

\begin{lemma}\label{lem:rel}
  The following relations hold for $x$ in $WH_b$, $d$ in $H_b$, and
  $s$ in $\CS$:
  \begin{enumerate}
  \item $t_xt_d= t_{xd}$ and $t_dt_x=t_{dx}$ \label{lem:rel1}
  \item $t_st_x= \begin{cases} t_{sx} &\text{if $\ell(sx)>\ell(x)$}\\
      \q t_{sx}+ \a \sum_{d\in X_s} t_{dx} &\text{if
        $\ell(sx)<\ell(x)$} \end{cases}$  \label{lem:rel2}
  \item $t_xt_s= \begin{cases} t_{xs} &\text{if $\ell(xs)>\ell(x)$}\\ 
      \q t_{xs}+ \a \sum_{d\in X_s} t_{xd} &\text{if
        $\ell(xs)<\ell(x)$} \label{lem:rel3}
    \end{cases}$
  \end{enumerate}
\end{lemma}

\begin{proof}
  The first statement follows easily from relation (\ref{rel1.5}) and
  the definitions.

  Using induction on $\ell(w)$ and (\ref{rel4}) it is easily seen that
  (\ref{lem:rel2}) holds when $x$ is replaced by $w$. Then
  (\ref{lem:rel2}) follows in general by writing $x=wd$ and using
  (\ref{lem:rel1}).
  
  Using induction on $\ell(w)$ and (\ref{rel4}) it is easily seen that
  (\ref{lem:rel3}) holds when $x$ is replaced by $w$. Then
  (\ref{lem:rel3}) follows in general by writing $x=dw$ and using
  (\ref{lem:rel1}).
\end{proof}

\begin{theorem}\label{generic}
  The algebra $\CH$ is free as an $A$-module with basis $\{\, t_x\mid
  x\in WH_b\,\}$.
\end{theorem}

\begin{proof}
  It follows from the last lemma that the span of $\{\, t_x\mid x\in
  WH_b\,\}$ is a two-sided ideal in $\CH$ containing the identity
  element so $\{\, t_x\mid x\in
  WH_b\,\}$ spans $\CH$. 
  
  To show that $\{\, t_x\mid x\in WH_b\,\}$ is linearly independent, we
  adapt Lusztig's presentation \cite[Proposition 3.3]{lusztig:hecke} of the
  argument sketched in \cite[Ch.~4, Ex.~23]{bourbaki:groupes} to the current
  situation.
  
  Let $E$ be a free $A$-module with basis $\{\, e_x\mid x\in
  WH_b\,\}$. For $s$ and $t$ in $\CS$ and $d$ in $H_b$, define
  endomorphisms $P_d$, $Q_d$, $P_s$, and $Q_t$ of $E$ by
  $A$-linearity and
  \begin{gather*}
    P_d(e_x)= e_{dx}, \qquad   Q_d(e_x)= e_{xd}\\
    P_s(e_x)= \begin{cases} e_{sx} &\ell(sx)>\ell(x)\\ \q e_{sx}+ \a
      \sum_{d'\in X_s} e_{d'x} &\ell(sx)<\ell(x), \end{cases}\\
    Q_t(e_x)= \begin{cases} e_{xt} &\ell(xt)>\ell(x)\\ \q e_{xt}+ \a
      \sum_{d'\in X_t} e_{xd'} &\ell(sx)<\ell(x). \end{cases}
  \end{gather*}
  
  We next show that $P_yQ_z=Q_zP_y$ for $y$ and $z$ in $\CS \cup H_b$.
  
  Clearly $P_dQ_{d'}= Q_{d'}P_d$, $P_sQ_d= Q_dP_s$, and $P_dQ_t=
  Q_tP_d$. The length function $\ell$ is constant on $H_b$ cosets and
  it follows that $P_sP_d= P_dP_s$ and $Q_dQ_t= Q_tQ_d$.
  
  It remains to show that $P_sQ_t= Q_tP_s$ for $s$ and $t$ in $\CS$. 

  Fix $x$ in $WH_b$. As in the proof of \cite[Proposition
  3.3]{lusztig:hecke} there are six cases. Suppose first that
  $\ell(sxt)=\ell(x)> \ell(sx)= \ell(xt)$. Then
  \[ 
  P_sQ_t(e_x)= \q e_{sxt}+ \a\q \sum_{d'\in X_t}e_{sxd'}+ \a^2
  \sum_{d\in X_s} \sum_{d'\in X_t} e_{dxd'}
  \]
  and
  \[ 
  Q_tP_s(e_x)= \q e_{sxt}+ \a\q \sum_{d\in X_s}e_{dxt}+ \a^2
  \sum_{d\in X_s} \sum_{d'\in X_t} e_{dxd'}.
  \]  
  Say $x=w\tilde d$ with $w$ in $W$ and $\tilde d$ in $H_b$. Then
  $\ell(swt)=\ell(w)> \ell(sw)= \ell(wt)$. By Deodhar's Property Z
  \cite[Theorem 1.1]{deodhar:bruhat} we have $wt=sw$ and so
  $s=wtw\inverse$. It follows from the definition that
  $X_s=wX_tw\inverse$. Moreover, it follows from the definition that
  $tX_t=X_tt$ and $X_t t\tilde d t= X_t \tilde d$. Therefore,
  $wt\tilde dX_t=wX_t\tilde dt$ and so $sxX_t=X_sxt$.  It follows that
  $P_sQ_t= Q_tP_s$.

  The other five cases are all easier. We omit the details.

  Now let $\CHt$ be the subalgebra of $\End_{A}(E)$ generated by $\{\,
  P_y\mid y\in \CS\cup H_b\,\}$. Consider the evaluation map,
  $\epsilon\colon \CHt\to E$, with $\epsilon(f)=f(e_1)$. We will show
  that $\epsilon$ is an isomorphism of $A$-modules. If $d$ is in $H_b$
  and $w=s_{i_1} \dotsm s_{i_p}$ is in $W$ with $\ell(w)=p$, then
  $P_{s_1} \dotsm P_{s_p} P_d(e_1)= e_{wd}$ and so $\epsilon$ is
  surjective. To show that $\epsilon$ is injective, suppose that $f$
  is in $\CHt$ and $f(e_1)=0$. If $d$ is in $H_b$ and $w=s_{i_1}
  \dotsm s_{i_p}$ is in $W$ with $\ell(w)=p$, then $f'=Q_dQ_{s_p}
  \dotsm Q_{s_1} (e_1)$ is in $\CHt$. By what we have shown above,
  $ff'=f'f$ and so
  \[
  0= f'f(e_1)= ff'(e_1)= f(e_{wd}).
  \]
  Since $wd$ is arbitrary in $WH_b$, it follows that $f=0$ and so
  $\epsilon$ is injective. Thus, $\epsilon$ is an isomorphism of
  $A$-modules.
  
  For $x$ in $WH_b$, let $f_x$ be the unique element in $\CHt$ with
  the property that $f_x(e_1)= e_x$. Then $\{\, f_x\mid x\in WH_b\,\}$
  is an $A$-basis of $\CHt$, $f_s=P_s$ for $s$ in $\CS$, and $f_d=P_d$
  for $d$ in $H_b$.  It is easily checked that the relations
  (\ref{rel1}) to (\ref{rel4}) are satisfied by the elements $\{\,
  f_s\mid s\in \CS\,\}\cup \{\, f_d\mid d\in H_b\,\}$ in $\CHt$. It
  follows that there is a homomorphism of $A$-algebras, $\phi\colon
  \CH\to \CHt$ with $\phi(t_y)=f_y$ for $y$ in $\CS \cup H_b$. Since
  $\{\, f_x\mid x\in WH_b\,\}$ is a basis of $\CHt$ it follows that
  $\{\,t_x\mid x\in WH_b\,\}$ is linearly independent. This completes
  the proof of the theorem.
\end{proof}

\begin{corollary}
  The element $t_1$ is the identity in $\CH$.
\end{corollary}

\begin{proof}
  By assumption, $t_1t_d=t_d=t_dt_1$ for $d$ in $H_b$ and
  $t_1t_w=t_wt_1$ for $w$ in $W$. Comparing coefficients in
  $t_d(t_1t_w)$ and $t_dt_w$ when both are expressed as linear
  combinations of $\{\, t_x\mid x\in WH_b\,\}$ shows that $t_1t_w=t_w$
  for $w$ in $W$. Therefore, $t_1t_x=t_x=t_xt_1$ for every $x$ in
  $WH_b$.
\end{proof}

Any function on $G$ that is constant on $(B,B)$-double cosets is
obviously constant on $(B_a, B_a)$-double cosets. Thus, if $e_B$ is
the centrally primitive idempotent in $\BBC GB$ corresponding to the
trivial representation of $B$, then $e_B\BBC Ge_B \subseteq e\BBC Ge
=\CHqa$. Consequently, one would expect that the Iwahori-Hecke algebra
of $W$ is a subalgebra of $\CH$. Taking into account the relation
$q-1=ab$ and that $\a$ is an indeterminate we show that this is indeed
the case.

The rule $d\mapsto t_d$ defines an $A$-algebra isomorphism between the
group algebra $AWH_b$ and the $A$-span of $\{\, t_d\mid d\in H_b\,\}$.
Define $e_1 =\sum_{d\in H_b} t_d$, so $t_de_1=e_1$ for every $d$ in
$H_b$ and $e_1^2= |H_b| e_1= b^n e_1$ (note that this is not the same
$e_1$ as in the proof above). It follows from (\ref{rel1.5}) that
$e_1$ is in the center of $\CH$ and so $\CH e_1$ is a two-sided ideal
in $\CH$.

\begin{corollary}
  The subalgebra $\CH e_1$ of $\CH$ is isomorphic to the Iwahori-Hecke
  algebra of $W$ with parameters $b^n\q$ and $b^{n+1}\,\a$.
\end{corollary}

\begin{proof}
  For $w$ in $W$ define $\tilde t_w=t_we_1$. It follows
  from Theorem \ref{generic} that $\{\, \tilde t_w\mid w\in W\,\}$ is
  a basis of $\CH e_1$. Clearly the elements $\tilde t_s$ for $s$ in
  $\CS$ satisfy the braid relations and
  \[
  (\tilde t_s)(\tilde t_s)= b^nt_s^2e_1=b^n \q t_1e_1+ b^n\a
  \sum_{d\in X_s} t_st_de_1= b^n\q \tilde t_1+ b^{n+1}\,\a \tilde t_s.
  \]
\end{proof}

In the rest of this section we record some results that follow more or less
immediately from Theorem \ref{generic}. For many of the constructions in the
rest of this paper it will be necessary to assume that $\q$ is invertible.
Also, some formulas become simpler if we rescale the basis elements $t_x$ of
$\CH$ by negative powers of $\bv$.  Thus, we let $A_\bv= \BBZ[\a,\bv,
\bv\inverse]$ denote the localization of $A$ at $\bv$ and define $\CHv=
A_\bv \otimes_A \CH$. For $x$ in $WH_b$, define $T_x= \bv^{-\ell(x)} \otimes
t_x$. Then $\{\, T_x\mid x\in WH_b\,\}$ is an $A_\bv$-basis of $\CHv$ and
the following relations hold for $x$ in $WH_b$, $d$ in $H_b$, and $s$ in
$\CS$:
\begin{gather}
  T_dT_{x}= T_{dx}\quad \text{and}\quad T_xT_d=T_{xd}.\label{vrel1}\\ 
  T_sT_x= \begin{cases} T_{sx} &\text{if $\ell(sx)>\ell(x)$}\\
    T_{sx}+ \a\bv\inverse\sum_{d\in X_s} T_{dx}&\text{if
      $\ell(sx)<\ell(x)$}. 
  \end{cases}\label{vrel2}\\
  T_xT_s= \begin{cases} T_{xs} &\text{if $\ell(xs)>\ell(x)$}\\
    T_{xs}+ \a\bv\inverse \sum_{d\in X_s} T_{xd} &\text{if
      $\ell(xs)<\ell(x)$}. \end{cases} \label{vrel3}
\end{gather}

The following lemma is a crucial ingredient needed to define the
$R$-polynomials and the partial order on $WH_b$.

\begin{lemma}\label{inverse}
  The elements $T_x$ for $x$ in $WH_b$ are units in $\CHv$. In
  particular
  \[
  T_s\inverse= T_s-\a\bv\inverse \sum_{d\in X_s} T_d
  \]
  for $s$ in $\CS$.
\end{lemma}

\begin{proof}
  Assume for a moment that $T_s$ is invertible for every $s$ in $\CS$.
  Then $T_w$ is invertible for every $w$ in $W$. It follows from
  (\ref{vrel1}) that $T_d$ is invertible for every $d$ in $H_b$ and
  thus $T_{wd}=T_wT_d$ is invertible for every $wd$ in $WH_b$. 
  
  To complete the proof, suppose $s$ is in $\CS$ and define $\Tsbar=
  T_s-\a\bv\inverse \sum_{d\in X_s} T_d$.  Then it follows from
  (\ref{vrel2}) and (\ref{vrel3}) that $\Tsbar$ is a left and right
  inverse for $T_s$.
\end{proof}

Define $\tau\colon \CHv\to A_\bv$ by $\tau(T_x)= \delta_{x,1}$ for $x$
in $WH_b$ and $A_\bv$-linearity, where $\delta_{x,y}$ is the Kronecker
delta.

\begin{proposition}\label{tau}
  If $x$ and $y$ are in $WH_b$, then $\tau(T_xT_y)=
  \delta_{x,y\inverse}$. Therefore, $\tau$ is a symmetrizing trace
  form on $\CHv$ and $\{\, T_{y\inverse}\mid y\in WH_b\,\}$ is the
  basis of $\CHv$ dual to the basis $\{\, T_x\mid x\in WH_b\,\}$.
\end{proposition}

\begin{proof}
  Assume that $\tau(T_xT_y)= \delta_{x, y\inverse}$. Then $\tau(hh')=
  \tau(h'h)$ for $h$ and $h'$ in $\CHv$ and the bilinear form
  $(h,h')\mapsto \tau(hh')$ is non-degenerate, so $\tau$ is a
  symmetrizing trace form on $\CHv$.
  
  We prove the formula $\tau(T_xT_y)= \delta_{x, y\inverse}$ using
  induction on $\ell(x)$.
 
  If $\ell(x)=0$, then $x=d$ is in $H_b$ and using (\ref{vrel1}) we
  have
   \[
   \tau(T_dT_y)= \tau(T_{dy})= \delta_{dy, 1}= \delta_{d, y\inverse}.
   \]
 
  Now suppose that $\ell(x)>0$. Write $x=dws$ where $d$ is in $H_b$,
  $w$ is in $W$, $s$ is in $\CS$, and $\ell(ws)>\ell(w)$. There are
  two cases: $\ell(sy)>\ell(y)$ and $\ell(sy)<\ell(y)$.
 
  First, suppose that $\ell(sy)>\ell(y)$. Since $\ell(xs)<\ell(x)$ and
  $\ell(y\inverse s)> \ell(y\inverse)$, it follows that $x\ne
  y\inverse$ and so $\delta_{x, y\inverse}=0$. Using (\ref{vrel1}),
  (\ref{vrel2}), and induction we have
  \[
  \tau(T_xT_y)= \tau(T_{dw}T_sT_y)= \tau(T_{dw}T_{sy}) = \delta_{dw,
    y\inverse s}.
  \]
  Since $x\ne y\inverse$ it follows that $dw\ne y\inverse s$.
  Therefore, $\tau(T_xT_y) =0$.
 
  Second, suppose that $\ell(sy)<\ell(y)$. Then using (\ref{vrel1}),
  (\ref{vrel2}), and induction we have
  \begin{align*}
    \tau(T_xT_y)&= \tau(T_{dw}T_sT_y)\\
    &= \tau \left(T_{dw} (T_{sy} +\a\bv\inverse \sum_{d_1\in X_s}
      T_{d_1y}) \right) \\
    &= \tau(T_{dw} T_{sy})+ \a\bv\inverse \sum_{d_1\in X_s} \tau(
    T_{dw}T_{d_1y}) \\
    &= \delta_{dw, y\inverse s} +\a\bv\inverse \sum_{d_1\in X_s}
    \q^{\ell(dw)} \delta_{ dw, (d_1y)\inverse}.
  \end{align*}
  If $dw=(d_1y)\inverse$ for some $d_1$ in $X_s$, then $\ell(y\inverse
  d_1\inverse s)= \ell(dws)> \ell(dw)= \ell( y\inverse d_1\inverse)$,
  so $\ell(y\inverse s)> \ell( y\inverse)$, a contradiction.
  Therefore, $dw\ne (d_1y)\inverse$ for $d_1$ in $X_s$ and so
  \[
  \tau(T_xT_y)= \delta_{dw, y\inverse s}= \delta_{x, y\inverse}
  \]
  as desired.
\end{proof}

Next, recall that a \emph{specialization} of $A_\bv$ is a ring
homomorphism from $A_\bv$ to a commutative ring $B$.  Given a
specialization, $\phi\colon A_\bv\to B$, we can consider the
$B$-algebra $\CHv\otimes _{A_\bv} B$.

Let $\phi_{1,0}\colon A_\bv\to \BBC$ be a ring homomorphism with
$\phi_1(\a)=0$ and $\phi_1(\bv)=\pm 1$.  Then clearly the
$\BBC$-algebra $\CH_{1,0}= \CHv\otimes_{A_\bv} \BBC$ is
isomorphic to the group algebra $\BBC WH_b$.

Let $\phi_{q,a}\colon A_\bv\to \BBC$ be a ring homomorphism with
$\phi_{q,a}(\a)=a$ and $\phi_{q,a}(\bv)=\pm \sqrt q$.  Then clearly the
$\BBC$-algebra $\CHv\otimes_{A_\bv} \BBC$ is isomorphic to the
algebra $\CHqa$ from \S2.

The arguments in \cite[\S68]{curtisreiner:methodsII} involving Tit's
Deformation Theorem, together with the fact that a group algebra is
always isomorphic to its opposite algebra, prove the following
theorem.

\begin{theorem}
  If $K$ is a field with characteristic zero that is a splitting field
  for $KWH_b$ and with the property that $\ind_{B_a}^G 1_{B_a}$ is a
  split $KG$-module, then the algebras $\End_{K G}(\ind_{B_a}^G
  1_{B_a})$ and $KWH_b$ are isomorphic.
\end{theorem}

Now suppose that $l$ is a prime that divides $q-1$ and does not divide
$n!$. Say $l^r|q-1$ and $l^{r+1} \not| q-1$. Set $b= l^r$. Then $H_b$ is a
Sylow $l$-subgroup of $G$. Let $(\CO, k, K)$ be a sufficiently large
$l$-modular system for $G$ and let $\CB$ denote the principal block of $\CO
G$. We assume that $q$ is a square in $K$ and we consider $\CB$ as a
two-sided ideal in $\CO G$. Let $\phi\colon A_\bv\to K$ be a ring
homomorphism with $\phi(\a)=a$ and $\phi(\bv)=\pm \sqrt q$.

\begin{corollary}
  The $\CO$-module $\ind_{B_a}^G 1_{B_a}$ is a biprojective $\CB
  \times (\CH\otimes_{A_\bv}K)$-bimodule and induces a Morita
  equivalence between the principal block, $\CB$, of $\CO G$ and the
  specialized algebra $\CH\otimes_{A_\bv} K$.
\end{corollary}

\begin{proof}
  This follows from the theorem using the argument in \cite[Chapter
  23]{cabanesenguehard:representations}.
\end{proof}

\section{$R$-polynomials and a Partial Order on $WH_b$}

In this section we follow the constructions of Kazhdan and Lusztig in
\cite{kazhdanlusztig:coxeter} and define a ``bar'' involution of $\CHv$ and
$R$-polynomials in $A_\bv$. We use the non-vanishing of the $R$-polynomials
to define a relation on $WH_b$ and we show that this relation is a partial
order on $WH_b$. This partial order reduces to the Bruhat order in case
$b=1$.

Define a ring endomorphism $\overline{\phantom{x}} \colon A_\bv\to A_\bv$ by
$\BBZ$-linearity, $\vbar= \bv\inverse$, and $\abar= -\a\bv^{-2}$.  Notice
that
\begin{itemize}
\item $\overline{\abar}= \overline{-\a\bv^{-2}}= \a\bv^2\bv^{-2}=\a$ so
  $\overline{\phantom{x}}$ is an involution, and
\item $\overline{\a\bv\inverse} =-\a\bv\inverse$.
\end{itemize}

Next, extend $\overline{\phantom{x}}$ to an endomorphism of $\CHv$ by
defining
\[
\overline{\sum_{x\in WH_b} \gamma_xT_x}= \sum_{x\in WH_b} \overline{
  \gamma_x} T_{x\inverse}\inverse = \sum_{x\in WH_b} \overline{ \gamma_x}
\overline{T_x}
\]
where $\overline{T_x}$ is defined to be $T_{x\inverse} \inverse$. Notice
that $\overline{T_d}= T_d$ for $d$ in $H_b$.

\begin{proposition}
  The endomorphism $\overline{\phantom{x}} \colon \CHv\to \CHv$ is a ring
  isomorphism of order two.
\end{proposition}

\begin{proof}
  Clearly $\overline{\phantom{x}}$ has order two and so is a
  bijection. Since $\CHv$ is generated by $\{\, T_d\mid d\in H_b\,\} \cup
  \{\, T_s\mid s\in \CS\,\}$, to show it is a ring homomorphism, it is
  enough to show that
  \[
  \overline{T_d T_x}= \overline{T_{d}} \,\overline {T_x}\quad
  \text{and}\quad \overline{T_s T_x}= \overline{T_s} \,\overline {T_x}
  \]
  for $x$ in $WH_b$, $d$ in $H_b$, and $s$ in $\CS$. These equations are
  easily shown to be true using (\ref{vrel1}) and (\ref{vrel2}).
\end{proof}

For $y$ in $WH_b$ we may express $\overline{T_y}$ in terms of the basis
$\{\, T_x\mid x\in WH_b\,\}$ of $\CHv$. For $x$ and $y$ in $WH_b$, define
elements $R^*_{x,y}$ in $A_\bv$ by
\[
\overline{T_y}= \sum_{x\in WH_b} \overline{R^*_{x,y}}\, T_x.
\]

Clearly $R^*_{1,1}=1$. If $d$ is in $H_b$, then $\overline{T_d}=T_d$ and so
$R^*_{x,d}=0$ if $x\ne d$ and $R^*_{d,d}=1$. It follows from Lemma
\ref{inverse} that for $s$ in $\CS$ we have $R^*_{x,s}=0$ unless $x$ is in
$\{\,s\,\} \cup X_s$, $R^*_{s,s}=1$, and $R^*_{d,s}=\a\bv\inverse$ for $d$
in $X_s$.

\begin{proposition}\label{recursion}
  Suppose $x$ and $y$ are in $WH_b$.
  \begin{enumerate}
  \item If $d$ is in $H_b$, then $R^*_{x, yd}= R^*_{xd\inverse,y}$ and
    $R^*_{x, dy}= R^*_{d\inverse x,y}$.
  \item If $s$ is in $\CS$ and $\ell(sy)<\ell(y)$, then
    \[
    R^*_{x,y}= \begin{cases}R^*_{sx,sy}&\text{if $\ell(sx)<\ell(x)$} \\
      R^*_{sx,sy}+ \a\bv\inverse \sum_{d\in X_s} R^*_{dx,sy}&\text{if
        $\ell(sx)>\ell(x)$}.
    \end{cases}
    \]
  \item If $s$ is in $\CS$ and $\ell(ys)<\ell(y)$, then
    \[
    R^*_{x,y}= \begin{cases}R^*_{xs,ys}&\text{if $\ell(xs)<\ell(x)$} \\
      R^*_{xs,ys}+ \a\bv\inverse \sum_{d\in X_s} R^*_{xd,ys}&\text{if
        $\ell(sx)>\ell(x)$}.
    \end{cases}
    \]
  \end{enumerate}
\end{proposition}

\begin{proof}
  To prove the first statement we compute $\overline{T_{yd}}$. Using
  (\ref{vrel1}) and the fact that $\overline{T_d}= T_d$ we have
  \[
  \overline{T_{yd}}= \overline{T_{y}} T_d = \sum_x
  \overline{R^*_{x,y}}T_xT_d = \sum_x \overline{R^*_{x,y}}T_{xd} = \sum_x
  \overline{R^*_{xd\inverse,y}}T_{x}.
  \]
  On the other hand, $\overline{T_{yd}}= \sum_x \overline{R^*_{x,yd}}\,
  T_x$. Comparing coefficients of $T_x$ we see that $R^*_{x, yd}=
  R^*_{xd\inverse,y}$.
  
  A similar argument shows that $R^*_{x, dy}= R^*_{d\inverse x,y}$.
    
  Next, by Lemma \ref{inverse}, (\ref{vrel1}), and (\ref{vrel2}) we have
  \begin{equation*}
    \overline{T_{y}}= \overline{T_{s}}\, \overline{T_{sy}} =
    \sum_{\substack{x\\ \ell(sx)<\ell(x)}} \overline{R^*_{sx,sy}} T_{x} +
    \sum_{\substack{x\\ \ell(sx)>\ell(x)}}\left(
      \overline{R^*_{sx,sy}}-\a\bv\inverse \sum_{d\in X_s}
      \overline{R^*_{d\inverse x,sy}} \right )T_{x}.    
  \end{equation*}
  The second statement of the proposition follows by comparing coefficients
  of $T_x$ and using that $\overline{\a\bv\inverse} =-\a\bv\inverse$ and
  $X_s=X_s\inverse$.
    
  The proof of the third statement is similar to the proof of the second
  statement and is omitted.
\end{proof}

\begin{corollary}\label{positive}
  For any $x$ and $y$ in $WH_b$, if $R^*_{x,y}\ne 0$, then $\ell(y)\geq
  \ell(x)$ and $R^*_{x,y}$ is a polynomial in $\a\bv\inverse$ with
  non-negative integer coefficients and degree at most $\ell(y)-\ell(x)$.
  In particular, $R^*_{x,y}=0$ if $\ell(x)\not\leq \ell(y)$. Moreover,
  $R^*_{y,y}=1$.
\end{corollary}

We will obtain a closed form for the polynomials $R^*_{x,y}$ in
\S\ref{subexp}.

\begin{proof}
  We prove the statements using induction on $\ell(y)$.
  
  If $\ell(y)=0$, then $y=d$ is in $H_b$ and we have seen that $R^*_{x,d}=0$
  unless $x=d$ and that $R^*_{d,d}=1$.
  
  Now suppose $\ell(y)>0$. Choose $s$ in $\CS$ with $\ell(sy)
  <\ell(y)$. Then by Proposition \ref{recursion} and induction,
  $R^*_{y,y}=R^*_{sy,sy}=1$.

  Now suppose $x$ is in $WH_b$ and $R^*_{x,y}\ne0$.
  
  If $\ell(sx)<\ell(x)$, then $R^*_{x,y}= R^*_{sx,sy}$ by Proposition
  \ref{recursion}. Since $\ell(y)- \ell(x)= \ell(sy)- \ell(sx)$, it follows
  by induction that $R^*_{x,y}$ is a polynomial in $\a\bv\inverse$ with
  non-negative integer coefficients and degree at most $\ell(y)-\ell(x)$.
  
  If $\ell(sx)>\ell(x)$, then $R^*_{x,y}= R^*_{sx,sy}+ \a\bv\inverse
  \sum_{d\in X_s} R^*_{dx,sy}$ by Proposition \ref{recursion}, and so either
  $R^*_{sx,sy}\ne0$ or $R^*_{dx,sy}\ne0$ for some $d$ in $X_s$.  Since
  $\ell(sy)-\ell(sx)= \ell(y)-\ell(x)-2$ and $\ell(sy)-\ell(dx)=
  \ell(y)-\ell(x)-1$ for any $d$ in $H_b$, it follows by induction that
  $R^*_{x,y}$ is a polynomial in $\a\bv\inverse$ with non-negative integer
  coefficients and degree at most $\ell(y)-\ell(x)$.
\end{proof}

For $x$ and $y$ in $WH_b$, define
\[
x\leq y\qquad \text{if}\qquad R^*_{x,y}\ne0.
\]
The main result in this section is that $\leq$ is a partial order on
$WH_b$.

Notice that when $b=1$, then $WH_b=W$. It is pointed out in
\cite{kazhdanlusztig:coxeter} that $R^*_{x,y}\ne0$ if and only if $x$ is
less than or equal to $y$ in the Bruhat-Chevalley order and so the relation
$\leq$ is the Bruhat-Chevalley order in this case.

We have seen that for $s$ in $\CS$ and $x$ in $WH_b$, $R^*_{x,s}\ne0$ if and
only if $x$ is in $X_s$, so $x\leq s$ if and only if $x$ is in $X_s$. In
particular, if $1$ is the identity in $WH_b$, then $1\leq s$ if and only if
$b$ is odd. It follows that in general, the restriction of $\leq$ to $W$ is
not the Bruhat-Chevalley order on $W$.

Suppose $w_1$ and $w_2$ are in $W$ and $d_1$ and $d_2$ are in $H_b$.  Then
by Proposition \ref{recursion}, $R^*_{w_1d_1, w_2d_2}=
R^*_{w_1d_1d_2\inverse, w_2}$. Define
\[
\Omega_{w_1, w_2}=\{\, d\in H_b\mid R^*_{w_1d,w_2}\ne0\,\}.
\]
Then for $w_1$, $w_2$ in $W$ and $d_1$, $d_2$ in $H_b$ we have
\[
\text{$w_1d_1\leq w_2d_2$ if and only if $d_1d_2\inverse \in
  \Omega_{w_1, w_2}$.}
\]
In order to understand the relation $\leq$ we need to describe the subsets
$\Omega_{w_1, w_2}$ of $H_b$.

\begin{lemma}\label{omegarecursion}
  Suppose $w_1$ and $w_2$ are in $W$, $s$ is in $\CS$, and
  $\ell(sw_2)<\ell(w_2)$. Set $t=w_1\inverse sw_1$. Then
  \[
  \Omega_{w_1, w_2}= \begin{cases} \Omega_{sw_1, sw_2}& \text{if
      $\ell(sw_1)<\ell(w_1)$} \\
    \Omega_{sw_1, sw_2} \cup \Omega_{w_1, sw_2}X_t & \text{if
      $\ell(sw_1)>\ell(w_1)$}. \end{cases}
  \]
\end{lemma}

\begin{proof}
  If $d$ is in $H_b$, then
  \[
  R^*_{w_1d,w_2}= \begin{cases}R^*_{sw_1d,sw_2}&\text{if
      $\ell(sw_1)<\ell(w_1)$} \\
    R^*_{sw_1d,sw_2}+ \a\bv\inverse \sum_{d'\in X_s}
    R^*_{d'w_1d,sw_2}&\text{if $\ell(sw_1)>\ell(w_1)$}.
    \end{cases}
  \]
  By definition $X_t=w_1\inverse X_s w_1$ so we may rewrite the last
  equation as
  \[
  R^*_{w_1d,w_2}= \begin{cases}R^*_{sw_1d,sw_2}&\text{if
      $\ell(sw_1)<\ell(w_1)$} \\
    R^*_{sw_1d,sw_2}+ \a\bv\inverse \sum_{d_1\in X_t}
    R^*_{w_1d_1d,sw_2}&\text{if $\ell(sw_1)>\ell(w_1)$}.
    \end{cases}
  \]
  It follows immediately that $\Omega_{w_1, w_2}= \Omega_{sw_1, sw_2}$ if
  $\ell(sw_1)<\ell(w_1)$.
  
  Suppose $\ell(sw_1)>\ell(w_1)$. Since $R^*_{x,y}$ is in $\BBN[ \a\bv
  \inverse]$, it follows that $R^*_{w_1d, w_2}\ne0$ if and only if
  $R^*_{sw_1d, sw_2}\ne0$ or there is a $d_1$ in $X_t$ with $R^*_{w_1d_1d,
    sw_2}\ne0$. Clearly $d_1d$ is in $\Omega_{w_1, sw_2}$ if and only if $d$
  is in $\Omega_{w_1, sw_2} d_1\inverse$. Since $X_t\inverse =X_t$ we have
  \[
  \Omega_{w_1, w_2}= \Omega_{sw_1, sw_2} \cup \left( \cup_{d_1\in X_t}
    \Omega_{w_1, sw_2} d_1\inverse \right)= \Omega_{sw_1, sw_2} \cup
  \Omega_{w_1, sw_2} X_t\inverse= \Omega_{sw_1, sw_2} \cup \Omega_{w_1,
    sw_2} X_t.
  \]
\end{proof}

In the next proposition, $1$ denotes the identity in $WH_b$ and $\leq_B$
denotes the Bruhat-Chevalley order on $W$.

\begin{proposition}\label{omega}
  Suppose $w_1$ and $w_2$ are in $W$. Then the following statements hold:
  \begin{enumerate}
  \item $\Omega_{w_1, w_2}\ne \emptyset$ if and only if $w_1\leq_B
    w_2$. \label{omega1}
  \item $\Omega_{w_1,w_1}= \{\, 1\,\}$. \label{omega2}
  \item If $w_1<_B w_2$, then $\Omega_{w_1, w_2}= X_{t_1} \dotsm X_{t_r}$
    where $t_1$, \dots, $t_r$ is any sequence of reflections in $\CT$ with
    $w_2=w_1t_1\dotsm t_r$ and $\ell(w_1t_1\dotsm t_i)= \ell(w_1)+i$ for
    $1\leq i\leq r$. \label{omega3}
  \end{enumerate}
\end{proposition}

\begin{proof}
  We will prove the proposition using induction on $\ell(w_2)$.
   
  If $\ell(w_2)=0$, then $w_2=1$. We have seen that for $x$ in $WH_b$,
  $R^*_{x,1}=0$ if $x\ne 1$ and $R^*_{1,1}=1$. Thus $\Omega_{w,1}
  =\emptyset$ for $w\ne1$ and $\Omega_{1,1}=\{\,1\,\}$ and so
  (\ref{omega1}), (\ref{omega2}), and (\ref{omega3}) hold in this case.
   
  Now suppose that $\ell(w_2)>0$ and fix $s$ in $\CS$ with $\ell(sw_2)<
  \ell(w_2)$. Then by Lemma \ref{omegarecursion} we have $\Omega_{w_2, w_2}=
  \Omega_{sw_2, sw_2}$ and so by induction, $\Omega_{w_2, w_2}=
  \{\,1\,\}$. This shows that (\ref{omega2}) holds.
   
  Suppose $w_1$ is in $W$ and $w_1\ne w_2$.
   
  We consider first the case when $\ell(sw_1)< \ell(w_1)$.
  
  By Lemma \ref{omegarecursion} we have $\Omega_{w_1, w_2}= \Omega_{sw_1,
    sw_2}$ and so by induction, $\Omega_{w_1, w_2} \ne \emptyset$ if and
  only if $sw_1<_Bsw_2$. By Deodhar's Property Z we have $sw_1\leq_B sw_2$
  if and only if $w_1\leq_B w_2$. Therefore, $\Omega_{w_1, w_2} \ne
  \emptyset$ if and only if $w_1\leq_B w_2$ and so (\ref{omega1}) holds.
   
  Suppose $w_1<_B w_2$ and $t_1$, \dots, $t_r$ is any sequence of
  reflections in $\CT$ with $w_2=w_1t_1\dotsm t_r$ and $\ell(w_1t_1\dotsm
  t_i)= \ell(w_1)+i$ for $1\leq i\leq r$. Define $v_0=w_1$ and
  $v_i=w_1t_1\dotsm t_i$ for $1\leq i\leq r$. Notice that $v_r=w_2$.
  
  There are two possibilities: either $\ell(sv_i) < \ell(v_i)$ for $1\leq
  i\leq r$, or there is an $i$ with $\ell(sv_i)> \ell(v_i)$ and
  $\ell(sv_{i+1})< \ell(v_{i+1})$.
   
  Suppose that $\ell(sv_i) < \ell(v_i)$ for $1\leq i\leq r$.  Then $sw_2=
  sw_1t_1\dotsm t_r$ and $\ell(swt_1\dotsm t_i)= \ell(sw_1)+i$ for $1\leq
  i\leq r$. Using Lemma \ref{omegarecursion} and induction we have
  \[
  \Omega_{w_1, w_2}= \Omega_{sw_1, sw_2}= X_{t_1} \dotsm X_{t_r}.
  \]
  
  Now suppose there is an $i$ with $\ell(sv_i)> \ell(v_i)$ and
  $\ell(sv_{i+1})< \ell(v_{i+1})$. Then by Deodhar's Property Z, $sv_i
  \leq_B v_{i+1}$. But $\ell( sv_i) = \ell( v_{i+1})$ and so $sv_i=
  v_{i+1}$. Hence $sw_1t_1\dotsm t_i= w_1t_1\dotsm t_it_{i+1}$ and
  $w_1\inverse sw_1= t_1\dotsm t_it_{i+1}t_i \dotsm t_1$.
  
  Fix $i$ so that $i$ is maximal with $\ell(sv_i)> \ell(v_i)$.  Then
  $1<i<r$, $v_{i+1}= sv_i$, and for $j>i$ we have $\ell(sv_j)< \ell(v_j)$.
  Since $v_{i+1}= sv_i$ we have
  \[
  sv_{i+2}= sv_it_{i+1}t_{i+2}= v_{i+1}t_{i+1}t_{i+2}= v_1t_{i+2}.
  \]
  Set
  \[
  t_1'=w_1\inverse sw_1,\quad \text{and} \quad t_j'= \begin{cases}
    t_{j-1}&\text{for $2\leq j\leq i+1$} \\ t_j& \text{for $i+2\leq j\leq
      r$.} \end{cases}
  \]
  Then $sw_2=sw_1t_1'\dotsm t_r'$ and $\ell(sw_1t_1'\dotsm t_j')=
  \ell(sw_1)+j$ for $1\leq j\leq r$. Using the induction hypothesis, the
  fact that $w_1\inverse sw_1= t_1\dotsm t_it_{i+1}t_i \dotsm t_1$, and
  Proposition \ref{Xt} we have
  \begin{align*}
    \Omega_{w_1, w_2} &= \Omega_{sw_1, sw_2}\\
    &= X_{w_1\inverse sw_1} X_{t_1} \dotsm X_{t_i} X_{t_{i+2}} \dotsm
    X_{t_r} \\
    &= X_{t_1} \dotsm X_{t_i} X_{t_1\dotsm t_it_{i+1}t_i \dotsm t_1}
    X_{t_{i+2}} \dotsm X_{t_r} \\
    &= X_{t_1} \dotsm X_{t_r} .
  \end{align*}
  
  We have shown that (\ref{omega3}) holds in both cases and thus have
  completed the proof of the proposition when $\ell(sw_1)< \ell(w_1)$.
  
  For the remainder of the proof we assume that $\ell(sw_1)>
  \ell(w_1)$. Then by Lemma \ref{omegarecursion} we have $\Omega_{w_1, w_2}=
  \Omega_{sw_1, sw_2} \cup \Omega_{w_1, sw_2}X_t$ where $t=w_1\inverse
  sw_1$. By Deodhar's Property Z the following three conditions are
  equivalent:
  \begin{itemize}
  \item $sw_1\leq_B w_2$,
  \item $w_1\leq_B w_2$, and 
  \item $w_1\leq_B sw_2$.
  \end{itemize}
  If $w_1\leq_B w_2$, then $w_1\leq_B sw_2$ and so by induction
  $\Omega_{w_1, sw_2} \ne \emptyset$. Therefore $\Omega_{w_1, w_2} \ne
  \emptyset$. Conversely, if $\Omega_{w_1, w_2} \ne \emptyset$, then either
  $\Omega_{sw_1, sw_2} \ne \emptyset$ or $\Omega_{w_1, sw_2} \ne
  \emptyset$. In the first case, it follows by induction that $sw_1\leq_B
  sw_2$. But then $sw_1\leq_B w_2$ and so $w_1\leq_B w_2$.  In the second
  case, it follows by induction that $w_1\leq_B sw_2$ and so again
  $w_1\leq_B w_2$. This proves (\ref{omega1}).
  
  Now suppose that $w_1<_Bw_2$ and $t_1$, \dots, $t_r$ is any sequence of
  reflections in $\CT$ with $w_2=w_1t_1\dotsm t_r$ and $\ell(w_1t_1\dotsm
  t_i)= \ell(w_1)+i$ for $1\leq i\leq r$. As above, define $v_0=w_1$ and
  $v_i=w_1t_1\dotsm t_i$ for $1\leq i\leq r$.  As above we have $v_r=w_2$.
  
  We consider the subset $\Omega_{w_1, sw_2}X_t$ of $H_b$. Choose $i$
  maximal with $\ell(sv_i) > \ell(v_i)$. Then we have seen that $w_1
  \inverse sw_1= t_1\dotsm t_it_{i+1} t_i \dotsm t_1$ and $sv_{i+2}= v_i
  t_{i+2}$. Set
  \[
  t_j'= \begin{cases} t_{j}&\text{for $1\leq j\leq i$} \\ t_{j+1}&
    \text{for $i+1\leq j\leq r-1$.} \end{cases}
  \]
  Then $sw_2=w_1t_1'\dotsm t_{r-1}'$ and $\ell(w_1t_1'\dotsm t_j')=
  \ell(w_1)+j$ for $1\leq j\leq r-1$.  Using the induction hypothesis
  and Proposition \ref{Xt} we have
  \begin{align*}
    \Omega_{w_1, sw_2}X_t &= X_{t_1} \dotsm X_{t_i} X_{t_{i+2}}
    \dotsm X_{t_{r}} X_t\\
    &=X_{t_1} \dotsm X_{t_i} X_{t_1\dotsm t_it_{i+1}t_i \dotsm t_1}
    X_{t_{i+2}} \dotsm X_{t_{r}} \\
    &=X_{t_1} \dotsm X_{t_{r}}.
  \end{align*}
  
  There are two cases, either $sw_1\not \leq_B sw_2$ or $sw_1\leq_B sw_2$.
  
  Suppose that $sw_1\not \leq_B sw_2$. Then $\Omega_{sw_1, sw_2} =\emptyset$
  and so
  \[
  \Omega_{w_1, w_2}= \Omega_{w_1, sw_2}X_t =X_{t_1} \dotsm X_{t_{r}}.
  \]
  
  Finally, suppose that $sw_1 \leq_B sw_2$. Then $\Omega_{sw_1, sw_2}\ne
  \emptyset$. We will show that $\Omega_{sw_1, sw_2} \subseteq \Omega_{w_1,
    sw_2} X_t$ and so
  \[
  \Omega_{w_1, w_2}= \Omega_{sw_1, sw_2} \cup \Omega_{w_1, sw_2}X_t =
  \Omega_{w_1, sw_2}X_t= X_{t_1} \dotsm X_{t_{r}}.
  \]
  
  Choose reflections $t_2'$, \dots, $t_{r-1}'$ so that $sw_2=sw_1t_2' \dotsm
  t_{r-1}'$ and $\ell(sw_1 t_2' \dotsm t_i')= \ell(sw_1)+i-1$ for $2\leq
  i\leq r-1$. Then by induction $\Omega_{sw_1, sw_2}= X_{t_2'} \dotsm
  X_{t_{r-1}'}$. Set $t_1'=t$. Then $sw_2=w_1t_1' \dotsm t_{r-1}'$ and
  $\ell(w_1 t_1' \dotsm t_i')= \ell(w_1)+i$ for $1\leq i\leq r-1$ and so by
  induction
  \[
  \Omega_{w_1, sw_2}X_t= X_t X_{t_2'} \dotsm X_{t_{r-1}'}X_t=
  X_t\Omega_{sw_1, sw_2}X_t.
  \]
  The subset $X_tX_t$ of $H_b$ is a subgroup and so $\Omega_{w_1, sw_2}X_t$
  contains $\Omega_{sw_1, sw_2}$ as claimed.
  
  It follows that (\ref{omega3}) holds when $\ell(sw_1)< \ell(w_1)$.  This
  completes the proof of the proposition.
\end{proof}

\begin{theorem}
  The relation $\leq$ on $WH_b$ is a partial order.
\end{theorem}

\begin{proof}
  We have seen that $R^*_{x,x}=1$ for all $x$ in $WH_b$ and so the relation
  is reflexive.
  
  Suppose $x=w_1d_1$, $y=w_2d_2$, and $z=w_3d_3$ are in $WH_b$ with $w_1$,
  $w_2$, $w_3$ in $W$ and $d_1$, $d_2$, $d_3$ in $H_b$.
  
  If $x\leq y$ and $y\leq x$, then $d_1d_2\inverse$ is in $\Omega_{w_1,
    w_2}$, so by Proposition \ref{omega} we have $w_1\leq_B w_2$. Similarly,
  $w_2\leq_B w_1$ and so $w_1=w_2$. But then $\Omega_{w_1, w_2}= \{\,1\,\}$
  and so $d_1=d_2$. Therefore $x=y$ and so the relation is anti-symmetric.
  
  If $x\leq y$ and $y\leq z$, then $d_1d_2\inverse$ is in $\Omega_{w_1,
    w_2}$ and $d_2d_3\inverse$ is in $\Omega_{w_2, w_3}$.  It follows from
  Proposition \ref{omega} that $w_1\leq_B w_2$ and $w_2\leq_B w_3$. Thus
  $w_1\leq_B w_2 \leq_B w_3$ and so using Proposition \ref{omega} again we
  see that $\Omega_{w_1, w_3}= \Omega_{w_1, w_2} \Omega_{w_2, w_3}$. Hence
  $d_1d_3\inverse = d_1d_2\inverse d_2d_3\inverse$ is in $\Omega_{w_1,
    w_3}$. Therefore $x\leq z$ and so the relation is transitive.
\end{proof}

We conclude this section with some properties of the partial order and some
examples.

\begin{proposition}\label{propd}
  Multiplication by $d$ is a poset automorphism of $(WH_b, \leq)$ for every
  $d$ in $H_b$.
\end{proposition}

\begin{proof}
  Clearly $x\mapsto xd$ and $x\mapsto dx$ are bijective mappings and it
  follows from Proposition \ref{recursion} that $x\leq y$ if and only if
  either $dx \leq dy$ or $xd \leq yd$.
\end{proof}

Define $H'=X_{s_1} \times \dotsm \times X_{s_{n-1}}$ and recall that
$X_0=\{\, h_1(\zeta^{ai}) \mid 0\leq i\leq b-1\,\}$. By Proposition
\ref{product}, the multiplication mapping from $X_0\times H'$ to $H_b$ is a
bijection.

\begin{proposition}\label{contain}
  If $b$ is odd, then $H'$ is a subgroup of $H_b$ that contains $X_{t_1}
  \dotsm X_{t_r}$ for every $t_1$, \dots, $t_r$ in $\CT$.
\end{proposition}

\begin{proof}
  Since $b$ is odd, each $X_t$ is a subgroup of $H_b$ and so $X_{s_1}\dotsm
  X_{s_{n-1}}$ is a subgroup of $H_b$. To prove the proposition it is enough
  to show that $X_t$ is contained in $X_{s_1}\dotsm X_{s_{n-1}}$ for every
  $t$ in $\CT$.

  Suppose $t$ is in $\CT$ and $t$ interchanges the $i$th and $j$th standard
  basis vector of $\BBF_q^n$ with $i<j$. Then $X_t=\{\, h_{i,j}(\alpha) \mid
  \alpha\in F_b\,\}$. We may assume that $j-1>1$.  The result follows since
  \[
  h_{i,j}(\alpha)= h_{i,i+1}(\alpha) h_{i+1,i+2}(\alpha) h_{1+2, i+3}(
  \alpha) \dotsm h_{j-1,j}(\alpha).
  \]
\end{proof}

\begin{proposition}
  If $b$ is odd, then $WH'$ is a normal subgroup of $WH_b$ and the cosets of
  $WH'$ are the connected components of the Hasse diagram of the poset
  $(WH_b, \leq)$.
\end{proposition}

\begin{proof}
  It is straightforward to check that $WH'$ is a normal subgroup of $WH_b$
  and it follows from Proposition \ref{product} that $X_0$ is a complete set
  of coset representatives of $WH'$ in $WH_b$.
  
  It follows from Proposition \ref{propd} that the posets $(WH', \leq)$ and
  $(WH'd_0, \leq)$ are isomorphic for every $d_0$ in $X_0$.

  Suppose $w_1d_1$ and $w_2d_2$ are in $WH'$ and $\alpha_1$ and $\alpha_2$
  are in $F_b$ with $w_1d_1h_1(\alpha_1) \leq w_2d_2h_1(\alpha_2)$. Then
  $d_1d_2\inverse h_1(\alpha_1 \alpha_2\inverse)$ is in $\Omega_{w_1,
    w_2}$. By the last proposition and Proposition \ref{omega} we have
  $\Omega_{w_1, w_2} \subseteq H'$ and so it follows from Proposition
  \ref{product} that $h_1(\alpha_1 \alpha_2\inverse)=1$. Therefore,
  $h_1(\alpha_1)= h_1(\alpha_2)$. This shows that if $x$ and $y$ are in
  $WH'$ and $d_0$ and $d_0'$ are in $X_0$, then no element of $WH'd_0$ is
  related to any element of $WH'd_0'$ if $d_0\ne d_0'$.
  
  Finally, suppose $w_1d_1$ and $w_2d_2$ are in $WH'$. We have seen that
  $\Omega_{1,w}$ is a subgroup of $H'$ for every $w$ in $W$. It follows that
  $d_1\leq w_1d_1$ and $d_2\leq w_2d_2$. Also, every simple reflection in
  $\CS$ occurs in any reduced expression for the longest element $w_0$, so
  $\Omega_{1, w_0}= H'$ and hence $d_1\leq w_0d_1d_2$ and $d_2\leq w_0
  d_1d_2$. This shows that the Hasse diagram of the poset $(WH', \leq)$ is
  connected.
\end{proof}

If $b$ is even, then it is not hard to show that the Hasse diagram of
$(WH_b, \leq)$ still has $b$ connected components, but they are somewhat
more complicated to describe. This is illustrated in the next example.

\begin{example}
  Suppose $n=3$ and $b=2$. Then the Hasse diagram of $WH_b$ has two
  connected components. The connected component containing the identity in
  $WH_2$ is given in Figure \ref{d:a1}. In this diagram we have denoted
  $h_i(-1)$ simply by $d_i$ for $i=1$, $2$, $3$. By Proposition \ref{propd},
  the other connected component of the Hasse diagram is obtained by
  multiplying on the left by $d_1$.

  \begin{figure}
    \centering
    \begin{sideways}
      \xymatrix{
        && s_1s_2s_1d_1 \ar@{-}[dddll] \ar@{-}[dddl] \ar@{-}[dddrrrr]
        \ar@{-}[dddrrrrr]& s_1s_2s_1d_1d_2d_3 \ar@{-}[dddlll] \ar@{-}[dddll]
        \ar@{-}[dddl] \ar@{-}[ddd]& s_1s_2s_1d_3 \ar@{-}[dddll]
        \ar@{-}[dddl] \ar@{-}[ddd] \ar@{-}[dddr]& s_1s_2s_1d_2 \ar@{-}[dddl]
        \ar@{-}[ddd] \ar@{-}[dddr]
        \ar@{-}[dddrr]&&\\
        &&&&&&&\\
        &&&&&&&\\
        s_2s_1d_1d_2 \ar@{-}[dddddd] \ar@{-}[ddddddrr] \ar@{-}[ddddddrrrr]
        \ar@{-}[ddddddrrrrrrr]& s_2s_1d_1d_3 \ar@{-}[dddddd]
        \ar@{-}[ddddddrr] \ar@{-}[ddddddrrrr] \ar@{-}[ddddddrrrrr]&
        s_1s_2d_2d_3 \ar@{-}[ddddddll] \ar@{-}[ddddddr] \ar@{-}[ddddddrrr]
        \ar@{-}[ddddddrrrrr]& s_1s_2d_1d_3 \ar@{-}[ddddddll]
        \ar@{-}[ddddddl] \ar@{-}[ddddddr] \ar@{-}[ddddddrrr]& s_2s_1
        \ar@{-}[ddddddllll] \ar@{-}[ddddddl] \ar@{-}[dddddd]
        \ar@{-}[ddddddrr]& s_2s_1d_2d_3 \ar@{-}[ddddddllll] \ar@{-}[ddddddlll]
        \ar@{-}[dddddd] \ar@{-}[ddddddrr]& s_1s_2d_1d_2
        \ar@{-}[ddddddllllll] \ar@{-}[ddddddllll] \ar@{-}[ddddddl]
        \ar@{-}[dddddd]& s_1s_2 \ar@{-}[ddddddllllll] \ar@{-}[ddddddllll]
        \ar@{-}[ddddddlll]
        \ar@{-}[dddddd]\\
        &&&&&&&\\
        &&&&&&&\\
        &&&&&&&\\
        &&&&&&& \\
        &&&&&&& \\
        s_2d_2 \ar@{-}[dddrr] \ar@{-}[dddrrr]& s_2d_3 \ar@{-}[dddr]
        \ar@{-}[dddrr]& s_1d_1d_2d_3 \ar@{-}[dddr] \ar@{-}[dddrr]& s_1d_3
        \ar@{-}[ddd] \ar@{-}[dddr]& s_2d_1 \ar@{-}[ddd] \ar@{-}[dddr]&
        s_2d_1d_2d_3 \ar@{-}[dddl] \ar@{-}[ddd]& s_1d_1 \ar@{-}[dddllll]
        \ar@{-}[dddl]& s_1d_2 \ar@{-}[dddlllll] \ar@{-}[dddll]\\
        &&&&&&& \\
        &&&&&&&\\
        && 1& d_2d_3& d_1d_3&d_1d_2&& }
    \end{sideways}
    \caption{Identity component of the Hasse diagram for
      $G(2,1,3)$} \label{d:a1} 
  \end{figure}
\end{example}

\section{The Kazhdan-Lusztig Basis and Kazhdan-Lusztig Polynomials}

Recall that $a$, $b$, and $q$ are related by the equation $ab=q-1$. In this
section we define a ring of scalars for $\CH$ so that the equation
$b\a=\bv^2-1$ can be solved for $\a$. Using this new ring of scalars we can
define a Kazhdan-Lusztig basis of $\CH$ and Kazhdan-Lusztig polynomials for
$WH_b$ following the construction in \cite[\S5]{lusztig:hecke}.

Let $\BBZ_b$ denote the localization of $\BBZ$ at $b$ and define $I$ to be
the principal ideal in $\BBZ_b[\a, \bv]_\bv$ generated by
$b\a-(\bv^2-1)$. Then $\a\bv\inverse \equiv \bv/b-\bv\inverse/b \mod I$. We
set $\Atilde= \BBZ_b[\a, \bv]_\bv /I$. The restriction of the natural
projection $\BBZ_b[\a, \bv]_\bv \to \Atilde$ to the subring
$\BBZ_b[\bv]_\bv$ is a ring isomorphism. In the following we identify
$\Atilde$ with $\BBZ_b[\bv]_\bv$.

Now define $\CHtilde = \Atilde \otimes_A \CH$ and again denote
$\bv^{-\ell(x)} \otimes t_x$ by $T_x$. Then $\{\, T_x\mid x\in WH_b\,\}$ is
an $\Atilde$-basis of $\CHtilde$ and the quadratic relation (\ref{vrel2})
becomes
\begin{gather}
  T_sT_x= \begin{cases} T_{sx} &\text{if $\ell(sx)>\ell(x)$}\\
    T_{sx}+ (\bv/b-\bv\inverse/b) \sum_{d\in X_s} T_{dx}&\text{if
      $\ell(sx)<\ell(x)$}
  \end{cases}
\end{gather}
for $x$ in $WH_b$ and $s$ in $\CS$. Also, $\overline {b\a-(\bv^2-1)}=
-\bv^{-2}(b\a-(\bv^2-1))$ and so the bar involution passes to the quotient
$\Atilde$ and to $\CHtilde$. Notice that polynomials $R^*_{x,y}$ defined in
\S4 are now polynomials in the quantity $\bv/b-\bv\inverse/b$.

Define
\[
\Atilde_{\leq0}= \BBZ_b[\bv\inverse], \quad \Atilde_{<0} = \bv\inverse
\BBZ_b [\bv\inverse], \quad \CHtilde_{\leq0}= \oplus_{x\in WH_b}
\Atilde_{\leq0} T_x, \quad \text{and} \quad \CHtilde_{<0}= \oplus_{x\in
  WH_b} \Atilde_{<0} T_x.
\]

With this notation, the proof of Theorem 5.2 in \cite{lusztig:hecke} applies
word-for-word to prove the following theorem.

\begin{theorem}\label{klbasis}
  Suppose $y$ is in $WH_b$. Then there is a unique element $C_y$ in
  $\CHtilde$ such that (1) $\overline{C_y}=C_y$ and (2) $C_y\equiv T_y \mod
  \CHtilde_{<0}$.
\end{theorem}

Clearly $\{\, C_y\mid y\in WH_b\,\}$ is an $\Atilde$-basis of $\CHtilde$. 

Define polynomials $P_{x,y}^*$ for $x$ and $y$ in $WH_b$ by
\[
C_y=\sum_{x\in WH_b} P_{x,y}^* T_x.
\]
Then as in \cite[\S5]{lusztig:hecke} we have:
\begin{itemize}
\item $P_{x,y}^*\ne0$ implies that $x\leq y$.
\item $P_{y,y}^*=1$.
\item For $x<y$, $P_{x,y}^*$ is a polynomial in $\bv\inverse$ with
  $P_{x,y}^*(0)=0$ and degree (in $\bv\inverse$) equal to $\ell(y)-\ell(x)$.
\item $P_{x,y}^*$ may be computed recursively once the $R$-polynomials are
  known using the equation
  \[
  \overline{P_{x,y}^*}- P_{x,y}^*= R^*_{x,y}+ \sum_{\substack{z\\ x<z<y}}
    R^*_{x,z} P_{z,y}^*.
    \]
\end{itemize}

\begin{corollary}
  For $y$ in $WH_b$ and $d$ in $H_b$ we have $C_{yd}=C_y T_d$ and $C_{dy}=
  T_d C_y$.
\end{corollary}

\begin{proof}
  Recall first that $\overline{T_d}=T_d$ and so $\overline{C_yT_d}=
  C_yT_d$. Next, $T_xT_d=T_{xd}$ and $x\leq yd$ if and only if $xd\inverse
  \leq y$, so we have
  \[
  C_{y}T_d= \sum_{\substack{x\\ x\leq y}} P_{x,y} T_{xd} =
  \sum_{\substack{z\\ z\leq yd}} P_{zd\inverse,y} T_z = T_{yd} +
  \sum_{\substack{z\\ z< yd}} P_{zd\inverse,y} T_z.
  \]
  Therefore, $C_xT_d$ is in $\CHtilde_{\leq0}$ and $C_yT_d-T_{yd}$ is in
  $\CHtilde_{<0}$. It follows from Theorem \ref{klbasis} that
  $C_{yd}=C_yT_d$.

  A similar argument shows that $C_{dy}= T_d C_y$.
\end{proof}

\begin{corollary}\label{klp2}
  For $x$ and $y$ in $WH_b$ and $d$ in $H_b$ we have
  $P^*_{x,yd}=P^*_{xd\inverse,y}$ and $P^*_{x,dy}=P^*_{d\inverse x,y}$.
\end{corollary}

For $x$ and $y$ in $WH_b$, define 
\[
R_{x,y}= v^{\ell(y)- \ell(x)} R_{x,y}^* \quad \text{and} \quad P_{x,y}=
v^{\ell(y)- \ell(x)} P_{x,y}^*.
\]
Then $\overline{T_y}= \sum_{x\in WH_b} v^{\ell(y)- \ell(x)}
\overline{R_{x,y}}\, T_x$ and $C_y= v^{-\ell(y)} \sum_{x\in WH_b} P_{x,y}
T_x$.

\begin{proposition}\label{p2}
  For $x$ and $y$ in $WH_b$, $R_{x,y}$ is a polynomial in $\bv^2$ with
  degree (in $\bv^2$) at most $\ell(y)-\ell(x)$.
\end{proposition}

\begin{proof}
  This follows immediately from the recursion formula in Proposition
  \ref{recursion}.
\end{proof}

\begin{corollary}
  For $x$ and $y$ in $WH_b$, $P_{x,y}$ is a polynomial in $\bv^2$.
\end{corollary}

\begin{proof}
  This follows immediately from the last proposition and the recursion
  formula for $P^*_{x,y}$.
\end{proof}

The $\Atilde$-basis $\{\, C_y\mid y\in WH_b\,\}$ is a \emph{Kazhdan-Lusztig
  basis} of $\CHtilde$ and the polynomials $P_{x,y}$ are
\emph{Kazhdan-Lusztig polynomials} for $WH_b$.

\begin{example}
  We compute the Kazhdan-Lusztig $R$-polynomials $R_{x,y}$ and the
  Kazhdan-Lusztig polynomials $P_{x,y}$ when $W$ is the symmetric group $S_3$. 

  Since $R_{x,y}= P_{x,y} =0$ unless $x\leq y$, $R_{x,x}=P_{x,x} =1$,
  $R_{x,yd}=R_{xd\inverse, y}$ and $P_{x,yd}=P_{xd\inverse, y}$, it is
  enough to compute $R_{x,w}$ and $P_{x,w}$ for $x$ in $WH_b$ and $w$ in $W$
  with $x<w$. Notice that $x<w$ implies $\ell(x)< \ell(w)$.

  If $x<w$ and $\ell(w)-\ell(x)=1$, then
  \[
  R_{x,w}(\bv)= \bv^2/b-1/b \quad \text{and} \quad P_{x,y}(\bv)= 1/b.
  \]

  If $x<w$ and $\ell(w)-\ell(x)=2$, then
  \[
  R_{x,w}(\bv)= \bv^4/b^2-2\bv^2/b^2+ 1/b^2 \quad \text{and} \quad
  P_{x,y}(\bv)= 1/b^2.
  \]

  If $x<w$ and $\ell(w)-\ell(x)=3$, then $w=s_1s_2s_1$ and $\ell(x)=0$,
  so $x=d$ is in $\Omega_{1,s_1s_2s_1}= X_{s_1}^2 X_{s_2}$. There are two
  cases depending on whether or not $d$ is in $s_1X_2s_1$.
  \begin{align*}
    R_{d,s_1s_2s_1}(\bv)&=
    \begin{cases}
      \frac 1{b^2} \bv^6- \frac {3-b}{b^2} \bv^4+ \frac {3-b}{b^2} \bv^2-
      \frac {1}{b^2} & d\in s_1X_2s_1 \\
      \frac 1{b^2} \bv^6- \frac {3}{b^2} \bv^4+ \frac {3}{b^2} \bv^2- \frac
      {1}{b^2} & d\in X_{s_1}^2X_{s_2} \setminus s_1X_{s_2}s_1
    \end{cases} \\
    P_{d,s_1s_2s_1}(\bv)&=
    \begin{cases}
      \frac {b-1}{b^2} \bv^2+ \frac {1}{b^2}& d\in s_1X_2s_1 \\
      \frac {-1}{b^2} \bv^2+ \frac {1}{b^2} & d\in X_{s_1}^2X_{s_2}
      \setminus s_1X_{s_2}s_1
    \end{cases}
  \end{align*}
  Notice that in general $P_{x,y}$ does not have non-negative coefficients
  and that the degree of $P_{x,y}$ depends on $b$. In particular, if $b=1$,
  so $WH_b=W$, then $d=1$ is in $s_1X_{s_2} s_1$ and $P_{1, s_1s_2s_1}=1$
  whereas if $b>1$ and odd, then $P_{1, s_1s_2s_1}$ has degree $2$. 
\end{example}

\section{Subexpressions and a Closed Formula for $R_{x,y}$} \label{subexp}

In this section we adapt results of Deodhar \cite{deodhar:geometric} to
describe the lower order ideals in $WH_b$ and to give a closed form for the
polynomials $R^*_{x,y}$.

It follows from Proposition \ref{recursion} that for $y$ and $z$ in $WH_b$
and $d$ in $H_b$ we have
\begin{gather*}
  y\leq zd \quad \text{if and only if}\quad yd\inverse \leq z \quad
  \text{and}\\
  y\leq dz \quad \text{if and only if}\quad d\inverse y \leq z.
\end{gather*}

For $z$ in $WH_b$, define $\CL^z=\{\, x\in WH_b\mid x\leq z\,\}$. Then if
$z=wd$ with $w$ in $W$ and $d$ in $H_b$ we have $\CL^z=
\CL^wd\inverse$. Thus, to analyze the lower order ideals $\CL^z$ it is
enough to consider the case when $z=w$ is in $W$.

Recall that a tuple of elements of $\CS$, say $(s_1, \dots, s_p)$, is said
to be \emph{reduced} if $\ell(s_1\dotsm s_p)=p$.

Fix a reduced tuple, $\bs=(s_1, \dots, s_p)$, of elements of $\CS$. A
\emph{distinguished subexpression of $\bs$} is a $p+1$-tuple of elements of
$WH_b$, say $\bx=(x_0, x_1, \dots, x_p)$, with the properties
\begin{itemize}
\item[(DS1)] $x_0=1$,
\item[(DS2)] $x_{j-1}\inverse x_j\in \{\, s_j\,\}\cup X_{s_j}$ for
$1\leq j\leq p$, and
\item[(DS3)] $x_{j}\leq x_{j-1}s_j$ for $1\leq j\leq p$.
\end{itemize}
Let $\CD_\bs$ denote the set of distinguished subexpressions of $\bs$ and
let $\pi_\bs\colon \CD_\bs\to WH_b$ be the projection on the last factor,
$\pi_\bs(\bx)= x_p$.

In the rest of this section we frequently argue using induction on the
number of elements in $\bs$. If $\bs=(s_1, \dots, s_p)$ with $p>1$, define
$\bs'=(s_1, \dots, s_{p-1})$.

\begin{theorem}\label{lowerorder}
  If $\bs=(s_1, \dots, s_p)$ is a reduced tuple of elements of $\CS$ and
  $w=s_1\dotsm s_p$, then
  \[
  \CL^w=\{\, \pi_\bs(\bx)\mid \bx \in \CD_\bs\,\}.
  \]
\end{theorem}

\begin{proof}
  We prove the result using induction on $p$. If $p=0$ then the result is
  clear and if $p=1$ then the result follows from Proposition \ref{omega}.
  
  Suppose $p>1$. Set $s=s_p$.
  
  Choose $y$ in $WH_b$ with $y\leq w$. 
  
  If $ys\leq ws$, then $0\ne R^*_{y,w}= R^*_{ys,ws}$, so $ys\leq ws$. By
  induction, there is an $\bx'=(x_0, \dots, x_{p-1})$ in $\CD_{\bs'}$ with
  $\pi_{\bs'}(\bx')= ys$. Define $x_p=y$ and $\bx=(x_0, \dots, x_{p-1},
  x_p)$. Then $x_{p-1}\inverse x_p=sy\inverse y=s$ is in $\{\,s_p\,\} \cup
  X_{s_p}$ and $x_p=y=yss=x_{p-1} s$. Thus, $\bx$ is in $\CD_\bs$,
  $\pi_\bs(\bx)=y$, and so $y$ is in $\pi(\CD_{\bs})$ in this case.
  
  If $ys\not \leq ws$, then we must have $\ell(ys)> \ell(y)$ and $R^*_{y,w}=
  \a\bv\inverse \sum_{d\in X_s} R^*_{yd, ws}$. Thus, there is a $d$ in $X_s$
  with $yd\leq ws$. By induction, there is an $\bx'=(x_0, \dots, x_{p-1})$
  in $\CD_{\bs'}$ with $\pi_{\bs'}(\bx')= yd$.  Define $x_p=y$ and
  $\bx=(x_0, \dots, x_{p-1}, x_p)$. Then $x_{p-1}\inverse x_p= d\inverse
  y\inverse y= d\inverse$ is in $\{\,s_p\,\} \cup X_{s_p}$. Say $y=vd_1$
  where $v$ is in $W$ and $d_1$ is in $H_b$.  Then $\ell(vs)>\ell(v)$
  because $\ell(ys)>\ell(y)$. Now $d\inverse d_1\inverse d_1\leq s$, so
  $d_1\leq d_1ds$, and so $vd_1\leq vssd_1ds=vd_1ds$. It follows that
  $x_p=y=vd_1\leq vd_1ds=yds=x_{p-1}s$. Thus, $\bx$ is in $\CD_\bs$,
  $\pi_\bs(\bx)=y$, and so $y$ is in $\pi_\bs(\CD_{\bs})$ in this case also.
  
  We have shown that $\CL^w \subseteq \{\, \pi(\bx)\mid \bx \in
  \CD_\bs\,\}$.
  
  Conversely, suppose $\bx$ is in $\CD_\bs$. Then clearly $(x_0, \dots,
  x_{p-1})$ is in $\CD_{\bs'}$. Also, $x_{p-1} \inverse x_p$ is in $\{\,
  s\,\}\cup X_s$, so $x_p=x_{p-1}s$ or $x_p=x_{p-1}d$ for some $d$ in
  $X_s$. We need to show that $x_p\leq w$.
  
  Suppose first that $x_p=x_{p-1}s$ and $\ell(x_ps)<\ell(x_p)$. Then
  $R^*_{x_p,w}= R^*_{x_ps, ws}= R^*_{x_{p-1}, ws}$. By induction, $x_{p-1}
  \leq ws$. Therefore, $R^*_{x_{p-1}, ws}\ne0$ and so $R^*_{x_p,w}\ne0$.
  Hence $x_p\leq w$ in this case.
  
  Next, suppose that $x_p=x_{p-1}s$ and $\ell(x_ps)>\ell(x_p)$. Then
  \[
  R^*_{x_p,w}= R^*_{x_ps,ws}+ \a\bv\inverse \sum_{d\in X_s} R^*_{x_pd, ws}=
  R^*_{x_{p-1}, ws}+ \a\bv\inverse \sum_{d\in X_s} R^*_{x_pd, ws}.
  \]
  Since $(x_0, \dots, x_{p-1})$ is in $\CD_{\bs'}$ it follows by induction
  that $x_{p-1} \leq ws$ and so $R^*_{x_{p-1}, ws}\ne0$.  Therefore, it
  follows from Corollary \ref{positive} that $R^*_{x_p,w}\ne0$ and so
  $x_p\leq w$ in this case.
  
  Finally, suppose that $x_p=x_{p-1}d$ where $d$ is in $X_s$. Then
  $x_p=x_{p-1}d\leq x_{p-1}s$. By induction $x_{p-1}\leq ws$ and so
  $R^*_{x_{p-1}, ws} \ne0$. By Proposition \ref{recursion} we have
  \[
  R^*_{x_{p-1}s,w}= \begin{cases}R^*_{x_{p-1},ws}&\text{if
      $\ell(x_{p-1}s)<\ell(x_{p-1})$} \\
    R^*_{x_{p-1},ws}+ \a\bv\inverse \sum_{d\in X_s}
    R^*_{x_{p-1}sd,ws}&\text{if $\ell(x_{p-1}s)>\ell(x_{p-1})$}.
  \end{cases}
  \]
  It follows from Corollary \ref{positive} that $R^*_{x_{p-1}s, w}\ne0$ and
  so $x_{p-1}s\leq w$. Therefore, $x_p\leq w$ in this case also.
\end{proof}

Suppose $\bs=(s_1, \dots, s_p)$ is a reduced tuple and $\bx =(x_1, \dots,
x_p)$ is in $\CD_\bs$. Define
\[
I(\bx)=\{\,j\mid 1\leq j\leq p\text{ and } x_{j-1}\inverse x_j \in X_{s_j}
\,\},\quad \text{and}\quad n(\bx)= |I(\bx)|.
\]
The rest of this section is devoted to the proof of the following theorem.

\begin{theorem}\label{closed}
  Suppose $\bs=(s_1, \dots, s_p)$ is a reduced tuple and $w=s_1\dotsm
  s_p$. Then
  \[
  R^*_{y,w}= \sum_{\bx\in \pi_\bs\inverse(y)} (\a\bv\inverse)^{n(\bx)}
  \]
  for $y\leq w$.
\end{theorem}

We will prove the theorem using induction on $p$ and the recursion formula
from Proposition \ref{recursion}.  The argument is essentially the same as
that in Deodhar \cite{deodhar:geometric}, suitably modified so as to make
sense in our context.

Set $w=s_1\dotsm s_p$ and $s=s_p$. In the next lemmas, $p>1$ and $y$ denotes
an element in $WH_b$ with $y\leq w$. We will analyze the subsets $\pi_\bs
\inverse(y)$ of $\CD_{\bs}$.

Recall that we have defined $\bs'=(s_1, \dots, s_{p-1})$. It follows
immediately from (DS1) to (DS3) that $\bx'= (x_0, \dots, x_{p-1})$ is in
$\CD_{\bs'}$. Let $\theta$ denote the projection from $\CD_\bs$ to
$\CD_{\bs'}$ given by $\theta( \bx)= \bx'$.

Clearly, if $\bx$ and $\widetilde {\bx}$ are in $\pi_\bs\inverse(y)$ with
$\theta(\bx)= \theta(\widetilde {\bx})$, then $\bx= \widetilde {\bx}$. Thus,
$\theta|_{\pi_\bs\inverse(y)}$ is injective.

Suppose $\bx=(x_0, \dots, x_p)$ is in $\pi_\bs\inverse(y)$ and $p>1$.  Then
$x_{p-1}\inverse y$ is in $\{s\}\cup X_s$ and so
\begin{itemize}
\item $p$ is not in $I(\bx)$ if and only if $x_{p-1}=ys$ and
\item $p$ is in $I(\bx)$ if and only if $x_{p-1}=yd$ for some $d$ in $X_s$.
\end{itemize}
Set
\[
\pi_\bs\inverse(y)_p^*=\{\, \bx\in \CD_\bs\mid \pi_\bs(\bx)=y\text{ and }
p\notin I(\bx)\,\}
\]
and
\[
\pi_\bs\inverse(y)_p=\{\, \bx\in \CD_\bs\mid \pi_\bs(\bx)=y\text{ and } p\in
I(\bx)\,\}
\]
so $\pi_\bs\inverse(y)$ is the disjoint union of $\pi_\bs \inverse(y)_p$ and
$\pi_\bs \inverse(y)_p^*$.

\begin{lemma}\label{lem:ys<y}
  If $\ell(ys)< \ell(y)$, then 
  \begin{enumerate}
  \item $ys\leq ws$, 
  \item $\pi_\bs \inverse(y)= \pi_\bs\inverse(y)_p^*$, and
  \item $\theta|_{ \pi_\bs\inverse(y)}\colon \pi_\bs\inverse(y) \to
    \pi_{\bs'}\inverse(ys)$ is a bijection.
  \end{enumerate}
\end{lemma}

\begin{proof}
  By Proposition \ref{recursion} we have $R^*_{y,w}= R^*_{ys,ws}$ and so
  $ys\leq ws$.

  Say $\bx=(x_0, \dots, x_p)$ is in $\pi_\bs\inverse(y)$. Then $x_p=y$,
  $x_{p-1}\inverse y$ is in $\{s\}\cup X_s$, and $y\leq x_{p-1}s$. Just
  suppose that $x_{p-1}\inverse y=d$ is in $X_s$.  Then $y\leq yd\inverse s$
  and so $\ell(y)\leq \ell(ys)$, a contradiction. Therefore, $p$ is not in
  $I(\bx)$. This shows that $\pi_\bs \inverse(y)= \pi_\bs\inverse(y)_p^*$.

  To show that $\theta\left( \pi_\bs\inverse(y) \right)
  =\pi_{\bs'}\inverse(ys)$, suppose $\bx'=(x_0, \dots, x_{p-1})$ is in
  $\CD_{\bs'}$ and $x_{p-1}=ys$. Set $x_p=y$ and $\bx=(x_0, \dots,
  x_p)$. Then $x_{p-1}\inverse x_p=sy\inverse y=s$ and $x_p=y=yss=x_{p-1}s$,
  so $\bx$ is in $\pi_\bs\inverse(y)$ and $\theta(\bx)=\bx'$.
\end{proof}

\begin{lemma}\label{lt}
  Suppose $y\leq w$ and $y=vd$ with $v$ in $W$ and $d$ in $X_s$. If
  $\ell(ys)>\ell(y)$, then $y\leq yd_1s$ for all $d_1$ in $X_s$.
\end{lemma}

\begin{proof}
  Since $\ell(ys)>\ell(y)$ it follows that $\ell(vs)>\ell(v)$. If $d_1$ is
  in $X_s$, then $d_1\inverse\leq s$ and so $1\leq ds$. Thus, $d\leq (s)(sd
  d_1s)$ and so $y=vd\leq (vs)(sdd_1s)= yd_1s$.
\end{proof}

\begin{lemma}\label{lem:ys>yN}
  If $\ell(ys)> \ell(y)$ and $ys\not\leq ws$, then
  \begin{enumerate}
  \item $\pi_\bs \inverse(y)= \pi_\bs\inverse(y)_p$ and
  \item $\theta|_{ \pi_\bs\inverse(y)}\colon \pi_\bs\inverse(y) \to
    \cup_{d\in X_s} \pi_{\bs'}\inverse(yd)$ is a bijection.
  \end{enumerate}
\end{lemma}

\begin{proof}
  Say $\bx=(x_0, \dots, x_p)$. Then $x_p=y$ and $x_{p-1}\inverse y$ is in
  $\{s\}\cup X_s$. Just suppose that $x_{p-1}\inverse y=s$, so
  $y=x_{p-1}s$. Now $\bx'=(x_0, \dots, x_{p-1})$ is in $\CD_{\bs'}$ and so
  by Theorem \ref{lowerorder}, $x_{p-1}\leq ws$. But $x_{p-1}=ys$ and by
  assumption $ys\not\leq ws$, a contradiction.  Therefore, $p$ is not in
  $I(\bx)$. This shows that $\pi_\bs \inverse(y)= \pi_\bs\inverse(y)_p$.
  
  To show that $\theta\left( \pi_\bs\inverse(y) \right) =\cup_{d\in X_s}
  \pi_{\bs'}\inverse(yd)$, suppose $d$ is in $X_s$, $\bx'=(x_0, \dots,
  x_{p-1})$ is in $\CD_{\bs'}$, and $x_{p-1}=yd$.  Set $x_p=y$ and
  $\bx=(x_0, \dots, x_p)$. Then $x_{p-1}\inverse x_p=d\inverse y\inverse
  y=d\inverse$ is in $X_s$. It follows from the last lemma that $y\leq yds$
  and so $\bx$ is in $\pi_\bs\inverse(y)$ and $\theta(\bx)=\bx'$.
\end{proof}

\begin{lemma}\label{lem:ys>yY}
  Suppose $\ell(ys)\leq \ell(y)$ and $ys\leq ws$, then
  \[
  \theta|_{ \pi_\bs\inverse(y)_p^*}\colon \pi_\bs\inverse(y)_p^* \to
  \pi_{\bs'}\inverse(ys)\quad \text{and}\quad \theta|_{
    \pi_\bs\inverse(y)_p}\colon \pi_\bs\inverse(y)_p \to \cup_{d\in
    X_s} \pi_{\bs'}\inverse(yd)
  \]
  are bijections.
\end{lemma}

\begin{proof}
  If $\bx=(x_0, \dots, x_p)$ is in $\pi_\bs\inverse(y)_p^*$, then $x_p=y$
  and $p\notin I(\bx)$, so $x_{p-1}=ys$. Thus $\theta(\bx)$ is in
  $\pi_{\bs'}\inverse (ys)$.
  
  To show that $\theta \left( \pi_\bs\inverse(y)_p^*\right) =
  \pi_{\bs'}\inverse(ys)$, suppose $\bx'=(x_0, \dots, x_{p-1})$ is in
  $\pi_{\bs'} \inverse(ys)$. Set $x_p=y$ and $\bx=(x_0, \dots, x_p)$.  Then
  as above, $x_{p-1}\inverse x_p=sy\inverse y=s$ and $x_p=y=yss=x_{p-1}s$,
  so $\bx$ is in $\pi_\bs\inverse(y)_p^*$ and $\theta(\bx)=\bx'$.
  
  If $\bx=(x_0, \dots, x_p)$ is in $\pi_\bs\inverse(y)_p$, then $x_p=y$ and
  $p\in I(\bx)$, so $x_{p-1}=yd_1$ for some $d_1$ in $X_s$.  Thus
  $\theta(\bx)$ is in $\cup_{d\in X_s} \pi_{\bs'}\inverse (yd)$.
  
  To show that $\theta \left( \pi_\bs\inverse(y)_p\right) = \cup_{d\in X_s}
  \pi_{\bs'}\inverse(yd)$, suppose $d$ is in $X_s$ and $\bx'=(x_0, \dots,
  x_{p-1})$ is in $ \pi_{\bs'}\inverse(yd)$. Set $x_p=y$ and $\bx=(x_0,
  \dots, x_p)$.  Then as above, $x_{p-1} \inverse x_p= d\inverse$ and it
  follows from Lemma \ref{lt} that $y\leq yds$, so $\bx$ is in
  $\pi_\bs\inverse(y)_p$.
\end{proof}

\begin{proof}[Proof of Theorem \ref{closed}]
  We prove the theorem using induction on $p$. If $p=1$, then the result
  follows from (\ref{vrel3}).
  
  Suppose $p>1$, set $w=s_1\dotsm s_p$ and $s=s_p$, and suppose $y$ is in
  $WH_b$ with $y\leq w$. There are three cases:
  \begin{enumerate}
  \item $\ell(ys)<\ell(y)$,
  \item $\ell(ys)>\ell(y)$ and $ys\not\leq ws$, and 
  \item $\ell(ys)>\ell(y)$ and $ys\leq ws$.
  \end{enumerate}
  
  Suppose first that $\ell(ys)<\ell(y)$. Then $ys\leq ws$ and
  $n(\theta(\bx))= n(\bx)$ for $\bx$ in $\pi_\bs\inverse(y)$. Using
  Proposition \ref{recursion}, induction, and Lemma \ref{lem:ys<y} we have
  \begin{align*}
    R^*_{y,w}&= R^*_{ys, ws} \\
    &= \sum_{\bx'\in \pi_{\bs'}\inverse (ys)} (\a\bv\inverse
    )^{n(\bx')}\\ 
    &= \sum_{\bx\in \pi_{\bs}\inverse (y)} (\a\bv\inverse)^{n(
      \theta(\bx))}\\ 
    &= \sum_{\bx\in \pi_{\bs}\inverse (y)} (\a\bv\inverse)^{n(\bx)}.
  \end{align*}
  
  Second, suppose that $\ell(ys)>\ell(y)$ and $ys\not\leq ws$. Then
  $n(\bx)=n(\theta(\bx))+1$ for $\bx$ in $\pi_\bs\inverse(y)$.  Using
  Proposition \ref{recursion}, induction, and Lemma \ref{lem:ys>yN} we have
  \begin{align*}
    R^*_{y,w}&= \a\bv\inverse \sum_{d\in X_s} R^*_{yd, ws} \\
    &= \a\bv\inverse \sum_{d\in X_s} \sum_{\bx'\in \pi_{\bs'}\inverse
      (yd)} (\a\bv\inverse )^{n(\bx')}\\
    &= \sum_{d\in X_s} \sum_{\bx'\in \pi_{\bs'}\inverse (yd)}
    (\a\bv\inverse )^{n(\bx')+1}\\
    &= \sum_{\bx\in \pi_{\bs}\inverse (y)} (\a\bv\inverse)^{n(
      \theta(\bx))+1}\\
    &= \sum_{\bx\in \pi_{\bs}\inverse (y)} (\a\bv\inverse)^{n(\bx)}.
  \end{align*}
  
  Finally, suppose that $\ell(ys)>\ell(y)$ and $ys\leq ws$. Then
  $n(\bx)=n(\theta(\bx))$ for $\bx$ in $\pi_\bs\inverse(y)_p^*$ and
  $n(\bx)=n(\theta(\bx))+1$ for $\bx$ in $\pi_\bs\inverse(y)_p$. Using
  Proposition \ref{recursion}, induction, and Lemma \ref{lem:ys>yY} we have
  \begin{align*}
    R^*_{y,w}&=R^*_{ys, ws}+ \a\bv\inverse \sum_{d\in X_s} R^*_{yd, ws} \\
    &= \sum_{\bx'\in \pi_{\bs'}\inverse(ys)} (\a\bv\inverse)^{n(\bx')} +
    \a\bv\inverse \sum_{d\in X_s} \sum_{\bx'\in \pi_{\bs'}\inverse
      (yd)} (\a\bv\inverse )^{n(\bx')}\\
    &= \sum_{\bx'\in \pi_{\bs'}\inverse(ys)} (\a\bv\inverse)^{n(\bx')} +
    \sum_{d\in X_s} \sum_{\bx'\in \pi_{\bs'}\inverse (yd)}
    (\a\bv\inverse )^{n(\bx')+1}\\
    &= \sum_{\bx\in \pi_{\bs}\inverse(ys)_p^*} (\a\bv\inverse )^{n(
      \theta(\bx))} + \sum_{\bx\in \pi_{\bs}\inverse (y)_p}
    (\a\bv\inverse )^{n(\theta(\bx))+1}\\
    &= \sum_{\bx\in \pi_{\bs}\inverse(ys)_p^*} (\a\bv\inverse )^{ n(\bx)} +
    \sum_{\bx\in \pi_{\bs}\inverse (y)_p} (\a\bv\inverse
    )^{n(\bx)}\\
    &= \sum_{\bx\in \pi_{\bs}\inverse (y)} (\a\bv\inverse)^{n(\bx)}.
  \end{align*}
  This completes the proof of the theorem.
\end{proof}


\bibliographystyle{amsplain}
\providecommand{\bysame}{\leavevmode\hbox to3em{\hrulefill}\thinspace}
\providecommand{\MR}{\relax\ifhmode\unskip\space\fi MR }
\providecommand{\MRhref}[2]{%
  \href{http://www.ams.org/mathscinet-getitem?mr=#1}{#2}
}
\providecommand{\href}[2]{#2}


\end{document}